\documentclass{amsart}
\usepackage{amssymb}
\usepackage{pictex}
\usepackage{afterpage}
\newtheorem{thm}{Theorem}[section]
\newtheorem{cor}[thm]{Corollary}
\newtheorem{lem}[thm]{Lemma}
\newtheorem{prop}[thm]{Proposition}
\theoremstyle{definition}
\newtheorem{defn}[thm]{Definition}
\newtheorem{rem}[thm]{Remark}  
\newtheorem{rems}[thm]{Remarks}
\numberwithin{equation}{section}

\newcommand{\R}{\mathbb R}

\newcommand{\To}{\longrightarrow}

\newcommand{\Ck}{\mathbb{K}}

\newcommand{\C}{\mathbb{C}}
\newcommand{\Z}{\mathbb{Z}}
\newcommand{\inv}{^{-1}}

\newcommand{\x}{\times}

\renewcommand{\k}{\mathbb K}
\newcommand{\w}{{\bf w}}
\renewcommand{\v}{{\bf v}}

\newcommand{\B}{{\mathcal B}}
\newcommand{\Hom}{{\operatorname{Hom}}}

\begin{document}
\title[Parametrizations of flag varieties]
{Parametrizations of flag varieties}
\author{R. J. Marsh}%
\address{Department of Mathematics and
Computer Science, University of Leicester, University Road, Leicester LE1 7RH}%
\email{R.Marsh@mcs.le.ac.uk}%
\author{K. Rietsch}%
\address{Department of Mathematics,
            King's College London,
            Strand, London
            WC2R 2LS}%
\email{rietsch@mth.kcl.ac.uk}%
\thanks{
The first named author was supported by a
University of Leicester Research Fund Grant and a Leverhulme Fellowship.
The second named author is supported by a Royal Society Dorothy
Hodgkin Research Fellowship.}%
\date{February 11 2004, revised March 19 2004}
\begin{abstract}
For the flag variety $G/B$ of a reductive algebraic group $G$ we
define and describe explicitly a certain (set-theoretical)
cross-section $\phi: G/B\to G$.  The definition of $\phi$ depends
only on a choice of reduced expression for the longest element
$w_0$ in the Weyl group $W$. It assigns to any $gB$ a
representative $g\in G$ together with a factorization into simple root
subgroups and simple reflections. The cross-section $\phi$ is
continuous along the components of Deodhar's decomposition of
$G/B$ \cite{Deo:Decomp}.
We introduce a generalization of the Chamber Ansatz of
\cite{BeZel:TotPos} and give formulas for the factors of
$g=\phi(gB)$.
These results are then applied to parametrize explicitly the
components of the totally nonnegative part of the flag variety
$(G/B)_{\ge 0}$ defined by Lusztig \cite{Lus:TotPos94}, giving a
new proof of Lusztig's conjectured cell decomposition of
$(G/B)_{\ge 0}$. We also give minimal sets of inequalities
describing these cells.
\end{abstract}
\maketitle

\section{Introduction}
Consider a simply connected reductive algebraic group $G$ over
$\C$ (or Chevalley group over $\k$) with opposite Borel subgroups
$B^+$ and $B^-$. So for example $G=SL_d(\C)$ with the subgroups of
upper- and lower-triangular matrices. The flag variety $G/B^+$ may
be embedded in the projective space of a sufficiently general
representation of $G$, say $V=V(\rho)$, by
 $$
G/B^+\hookrightarrow \mathbb P(V)~: \ gB^+\to \left<g\cdot
\xi\right>_\C,
 $$
where $\xi$ is a highest weight vector. Then to any element
$gB^+$ we may associate the highest and lowest extremal weights,
$v\rho$ and $w\rho$, such that $g\cdot \xi$ has nonzero component
in the corresponding weight space. Equivalently, the Weyl group
elements $v$ and $w$ determine the intersection of opposed Bruhat
cells
 $$
 B^- v B^+/B^+\cap B^+ w B^+/B^+
 $$
in which $gB^+$ lies. Now fix a reduced expression
$w_0=s_{i_1}\cdots s_{i_N}$ for the longest element of the Weyl
group. Following V. Deodhar \cite{Deo:Decomp}, there is a
finer datum that can be associated to $gB^+$. The element $gB^+$
can be successively reduced, compatibly with this reduced
expression, to give a sequence $ \left(B^+, g_{(1)}B^+,\dotsc,
g_{(n-1)}B^+,gB^+\right )$ in the flag variety, or a sequence of
intermediate lines
 $$
L_0=\left<\xi\right>\, ,\, L_1=\left<g_{(1)}\cdot \xi\right>\, ,
\, \dotsc\, ,\,L_{n-1}=\left<g_{(n-1)}\cdot \xi\right>\, ,\,
L_{n}=\left<g\cdot \xi\right>
 $$
in $V(\rho)$. For example, if we write $gB^+$ as $bw B^+$ for $b\in
B^+$, then $L_{n-1}$ is the line $\left<bws_i\cdot \xi\right>$,
where $s_i$ is the right-most simple reflection in the reduced
expression for $w_0$ such that $ws_i<w$ (see
Section~\ref{s:DeoComps}). Given all the intermediate lines $L_k$,
the further data associated to $gB^+$ is now the collection
$(v_{(1)},\dotsc, v_{(n)})$ of Weyl group elements such that
$v_{(k)}\rho$ is the {\it highest} extremal weight for which
$g_{(k)}\cdot \xi$ has non-zero weight space component. The set
of $gB^+$ in $B^+w B^+/B^+$ with fixed $(v_{(1)},\dotsc,
v_{(n)})$ is called a Deodhar component of the flag variety.

Consider the special case where the element $gB^+$ from above has
$v=1$. Then $gB^+=uB^+$ for some unipotent $u\in B^-$. If also
$v_{(i)}=1$ for all $i$, then $u$ may be factorized into negative
simple root subgroups as $u=y_{j_1}(t_1)\cdots y_{j_n}(t_n)$ for
some nonzero parameters $t_i\in\C$ (where $s_{j_1}\cdots
s_{j_n}$ is a reduced expression for $w$ governing the
construction of the intermediate lines $L_i$). If we write
$uB^+=zw B^+$ for some unipotent
$z\in B^+$, then A.~Berenstein and A.~Zelevinsky's Chamber Ansatz
\cite{BeZel:TotPos} gives formulas for the $t_i$ in terms of
minors of $z$.

In this paper we generalize the above result by describing
factorizations, and hence parametrizations, for a general Deodhar
component and giving formulas for the parameters
(Proposition~\ref{p:parameterization} and Theorem~\ref{t:ansatz}).
Our formulas for the nonzero parameters, analogous to the $t_k$ above, are
obtained by a direct generalization of the Chamber Ansatz.
However a general Deodhar component also has
another type of parameter 
which runs through $\k$. The formulas for these involve the
generalized Chamber Ansatz along with a correction term.

The Chamber Ansatz used in the formulas for the parameters depends
on the Deodhar component in which an element $zwB^+$ lies. Therefore we
also give a simple algorithm to determine this component
(Section~\ref{s:algorithm}). The algorithm in a sense `generates'
the chambers in the Chamber Ansatz for $zwB^+$ recursively. We
illustrate how this works with a very explicit type $A$ example
in Section~\ref{s:componentexample}.

In Section~\ref{s:totpos} we set $\k=\R$ and use these results to
examine the totally nonnegative part $(G/B^+)_{\ge 0}$ of the
flag variety. This is the closure in $G/B^+$ of the set
$\{y_{i_1}(t_1)\cdots y_{i_N}(t_N)B^+\, |\, t_i\in\R_{>0}\}$. We
explicitly describe the intersection of $(G/B^+)_{\ge 0}$
with each of the sets
$\mathcal R_{v,w}=B^-v B^+/B^+\cap B^+w B^+/B^+$. Namely in
$\mathcal R_{v,w}$ there is a unique open dense Deodhar component
which is isomorphic to $(\R^*)^{\ell(w)-\ell(v)}$.
And the totally nonnegative part $\mathcal R^{>0}_{v,w}$
of $\mathcal R_{v,w}$ is shown to be the subset of the above
Deodhar component where all of the parameters are positive.

This in particular reproves a result of the second author
conjectured by G.~Lusztig, that $\mathcal R^{>0}_{v,w}$ is a
semi-algebraic cell. However the new proof presented in this paper
gives for the first time explicit parametrizations of these totally
nonnegative parts (depending on a choice of reduced expression of
$w$).
And it has the advantage of being independent of the theory of canonical
bases, which was required in the previous proof.
Moreover
the parameters of $\mathcal R_{v,w}^{>0}$ can all be computed by the
generalized Chamber Ansatz (without correction term).

Finally in Section~11 we give an efficient description for
$\mathcal R_{v,w}^{>0}$ in terms of minor inequalities,
generalizing a result of Berenstein and Zelevinsky from the $v=1$
case. For any choice of reduced expression for $w$ we obtain a
set of $\ell(w)-\ell(v)$ inequalities. This set of inequalities
and $\ell(v)$ minor equalities, that can also be given
explicitly, describe $\mathcal R_{v,w}^{>0}$ as a semi-algebraic
subset of the (real) Bruhat cell $B^+ w B^+/B^+$.

\begin{rem} The case of intersections of opposite Bruhat cells
$\mathcal R_{v,w}$ in the flag
variety which we treat in this paper is not to be confused with
intersections of opposite Bruhat double cosets
 $$G^{v,w}=B^+w B^+\cap B^-vB^-$$
in the group. These other intersections
were studied by Fomin and Zelevinsky \cite{FoZe:DoubleCells}, who
obtained a different generalization of the Chamber Ansatz in that
setting (using it also to give parametrizations and
minimal sets of inequalities for their corresponding totally
positive parts in the group).

Our study of parametrizations in flag varieties
compatible with the $\mathcal R_{v,w}$ is
substantially different from the problems
in the group considered in
\cite{FoZe:DoubleCells}, for example
already where total positivity is concerned.
An immediate and obvious difference between total positivity questions
in the two cases lies in the fact that $(G/B^+)_{\ge 0}$
is the {\it closure} of the image in $G/B^+$ of the totally nonnegative part
$G_{\ge 0}$ of the group, and is actually generally larger than this image.
So it is clear that the totally positive cells
$\mathcal R_{v,w}^{>0}$ in $(G/B^+)_{\ge 0}$ cannot all
come from totally positive cells in $G_{\ge 0}$.
In fact, the cell decomposition of $G_{\ge 0}$, which is
studied in detail in \cite{FoZe:DoubleCells},
was first obtained by Lusztig in \cite{Lus:TotPos94} where
the analogous problem for flag varieties was formulated only
as a conjecture. One has to depart significantly from
the study of total positivity in the group in order to
study total positivity in the flag variety.

The overlap between the two parametrization problems, ours and the one
from \cite{FoZe:DoubleCells}, is precisely the joint special case covered in
\cite{BeZel:TotPos}. In that case
one has $G^{1,w}\cong \mathcal R_{1,w}\x T$,
where the maximal torus factor $T$ is irrelevant for the
parametrization problem. Otherwise, unless $v=1$ or
symmetrically $w=w_0$, the varieties $G^{v,w}$ have no
sensible counterpart in the flag variety.
Moreover both \cite{BeZel:TotPos}
and \cite{FoZe:DoubleCells} parametrize and give
formulas only for an open dense subset of the varieties they study.
So Theorem~\ref{t:ansatz} already
adds to these results in the joint
special case,
since it determines parameters for {\it any} element
in $\mathcal R_{1,w}$.

It is an interesting open problem to extend our
results from $\mathcal R_{1,w}$, and hence $G^{1,w}$, also to
the remaining varieties $G^{v,w}$ in the group. That is,
similarly to find a way to parametrize every
element of the group $G$. This should involve finding
appropriate stratifications of
the $G^{v,w}$ for arbitrary $v,w$ (the open strata
being the ones already understood by
\cite{FoZe:DoubleCells}), and then extending the
Chamber Ansatz from \cite{FoZe:DoubleCells} to all
the remaining strata.
\end{rem}

\section{Notation and basic definitions}

Let $\k$ be a field. Let $G_\k$ be a split, connected, simply
connected, semisimple algebraic $\k$-group (or Chevalley group
over $\k$). See \cite{Jantzen:AlgGroupBook}~Section II.1 or any of
\cite{Borel:AlgGroupBook}, \cite{Springer:AlgGroupBook},
\cite{Steinberg:ChevGroupBook}. Fix a $\k$-split maximal torus
$T_\k$. We write $\k^*$ for the multiplicative group $G_m(\k)$ and
$\k$ for the additive group $G_a(\k)$. As we will always be
concerned with the $\k$-valued points we will write $G$ for
$G(\k)$ and $T$ for $T(\k)$, and so forth. In later sections we
will take $\k$ to be $\R$.

Let $X(T)=\Hom(T,\k^*)$ and $R\subset X(T)$ the set of roots.
Choose a system of positive roots $R^+$. We denote by $B^+$ the
Borel subgroup corresponding to $R^+$, and by $U^+$ its unipotent
radical. We also have the opposite Borel $B^-$ such that $B^+\cap
B^-=T$, and its unipotent radical $U^-$.

Denote the set of simple roots by
\begin{equation*}
 \Pi=\{\alpha_i\ |\ i\in I\}\subset R^+\subset R\subset X(T).
\end{equation*}
For every $\alpha_i\in\Pi$ there is an associated homomorphism
\begin{equation*}
\varphi_i: SL_2 \to G.
\end{equation*}
Consider the $1$-parameter subgroups in $G$ (landing in $U^+$,
$U^-$ and $T$ respectively) defined by
\begin{equation*}
x_i(m)=\varphi_i\begin{pmatrix}1&m\\ 0&1
\end{pmatrix},\quad
y_i(m)=\varphi_i\begin{pmatrix}1&0\\ m&1
\end{pmatrix},\quad
\alpha_i^\vee(t)=\varphi_i\begin{pmatrix}t&0\\ 0&t\inv
\end{pmatrix},
\end{equation*}
where $m\in\k, \ t\in\k^*$, and $i\in I$.  The datum $(T, B^+,
B^-, x_i, y_i; i\in I)$ for $G$ is called a {\it pinning} in
\cite{Lus:TotPos94}. The standard pinning for $SL_d$ consists of
the diagonal, upper-triangular and lower-triangular matrices,
along with the simple root subgroups $x_i(m)=I_d+m E_{i,i+1}$ and
$y_i(m)= I_d + m E_{i+1,i}$, where $I_d$ is the identity matrix,
and $E_{i,j}$ has a $1$ in position $(i,j)$ and zeroes elsewhere.

Next consider the cocharacter lattice $Y(T)=\Hom(\k^*,T)$. It is
dually paired with $X(T)$ in the standard way by $<\, ,\,
>:X(T)\x Y(T)\to\Hom(\k^*,\k^*)\cong\Z$. The
$\alpha_i^\vee$ viewed as elements of $Y(T)$ are the simple
coroots, and the Cartan matrix $A=(a_{ij})\in\Z^{I\x I}$ is given
by $a_{ij}=<\alpha_j,\alpha_i^\vee>$.

That $G$ is simply connected means the $\alpha_i^\vee$ freely generate
$Y(T)\cong\Z^I$. The dual basis in $X(T)$ is the set of
fundamental weights $\{\omega_i\ |\ i\in I\}$. Let $X(T)_+$ be the
set of dominant weights and $\rho=\sum_{i\in I}\omega_i\in
X(T)_+$. For a dominant weight $\lambda$ let $V(\lambda)$ denote
the Weyl module with highest weight $\lambda$, see
\cite{Jantzen:AlgGroupBook}~II~2.13. In characteristic $0$ this is
just the irreducible representation with highest weight $\lambda$.

The Weyl group $W=N_G(T)/T$ acts on $X(T)$ permuting the roots
$R$. We denote the action of $w\in W$ on $\alpha\in X(T)$ by
$w\alpha$. The simple reflections $s_i\in W$ are given explicitly
by $s_i:=\dot s_i T$, where
\begin{equation*}\label{e:si}
\dot s_i:=\varphi_i\begin{pmatrix}0&-1\\ 1&0\end{pmatrix},
\end{equation*}
and any $w\in W$ can be expressed as a product $w=s_{i_1}\cdots
s_{i_m}$ with a minimal number of factors $m=\ell(w)$. We set
 $$
 \dot w=\dot s_{i_1}\dot s_{i_2}\cdot \dotsc\cdot \dot s_{i_m}
 $$
to get a representative of $w$ in $N_G(T)$. It is well known that
this product is independent of the choice of reduced expression
$s_{i_1}\cdots s_{i_m}$ for $w$. Let $<$ denote the Bruhat order
on $W$. The unique maximal element of $W$ is denoted $w_0$.

We note for future reference the following identity
(\cite{Jantzen:AlgGroupBook}~II~1.3)
\begin{equation}\label{e:siIdentity}
\alpha_i^\vee(t\inv)\dot s_i=x_i(-t\inv)y_i(t)x_i(-t\inv),
\end{equation}
which can be checked in $SL_2(\k)$.

Finally, for every root we introduce the corresponding root
subgroup. Let $U^+_{\alpha_i}$ be the simple root subgroup in $G$
given explicitly by $\{x_{i}(t)\ |\ t\in \k\}$. For an arbitrary
root $\alpha$ there is a $w\in W$ and simple root $\alpha_i$ such
that $\alpha=w\alpha_i$. Then the one-dimensional subgroup
corresponding to $\alpha$ may be defined as $\dot w
U^+_{\alpha_i}\dot w\inv$. If $\alpha\in R^+$ this subgroup lies
in $U^+$ and we write $U^+_\alpha=\dot w U^+_{\alpha_i}\dot
w\inv$. Otherwise the subgroup is called  $U^-_\alpha=\dot w
U^+_{\alpha_i}\dot w\inv$ and lies in $U^-$.

\section{Subexpressions of reduced expressions}
Consider a reduced expression in $W$, say $s_3 s_2 s_1 s_3 s_2
s_3$ in type $A_3$. Informally, a subexpression is what is
obtained by choosing some of the factors. So for example choosing
the underlined factors in
 \begin{equation}\label{e:subexpr}
 \underline{s_3}\, \underline{s_2} \, s_1\, \underline{s_3
 }\, \underline{s_2}\,  s_3
 \end{equation}
gives a subexpression for $s_2 s_3$ in the word $s_3 s_2 s_1 s_3
s_2 s_3$.

It will be useful to represent expressions, like $s_3 s_2 s_1 s_3
s_2 s_3$ or its subexpression $s_3 s_2\, 1\, s_3 s_2\, 1$, by
their sequences of partial products
\begin{equation*}
\begin{array}{llllllllll}
&(&1,\, &s_3,\, &s_3 s_2,\, &s_3 s_2 s_1,\, &s_3 s_2 s_1 s_3, \,
&s_3 s_2 s_1 s_3
s_2,\, & s_3 s_2 s_1 s_3 s_2 s_3 &),\\
 &(&1,\, &s_3,\, &s_3 s_2,\, &s_3 s_2 ,\, &s_3 s_2 s_3 ,\,  &s_2
 s_3,\, &s_2 s_3&).
\end{array}
\end{equation*}
We formalize this below.
\begin{defn}
\label{d:expressions} Let us define an {\it expression} for $w\in
W$ to be a sequence
 $$\w=\left (w_{(0)},w_{(1)}, w_{(2)},\dotsc, w_{(n)}\right )$$
in $W$, such that $w_{(0)}=1$, $w_{(n)}=w$ and
 $$
w_{(j)}=\begin{cases}\ w_{(j-1)},&\text{or}\\
          \ w_{(j-1)}s_i,&\text{for some simple reflection $s_i$}
        \end{cases}
 $$
for $j=1,\dotsc,n$. The expression $\w$ may equivalently be
specified by its {\it sequence of factors},
 $$
 (w_{(1)},w_{(1)}\inv
w_{(2)}^{\ },\dotsc, w_{(n-1)}\inv w_{(n)}^{\ }),
 $$
which has entries in $\{s_i \ |\ i\in I \}\cup \{1\}$.
\end{defn}

\begin{defn}\label{d:Js}
For an expression $\w=(w_{(0)},w_{(1)},\dotsc, w_{(n)})$ define
\begin{align*}
J^+_\w &=\{k\in\{1,\dotsc,n\}\ |\  w_{(k-1)}<w_{(k)}\},\\
J^{\circ}_\w\, &=\{k\in\{1,\dotsc,n\}\ |\  w_{(k-1)}=w_{(k)}\},\\
J^-_\w &=\{k\in\{1,\dotsc,n\}\ |\  w_{(k)}<w_{(k-1)}\}.
\end{align*}
An expression $\w=(w_{(0)},w_{(1)},\dotsc,w_{(n)})$ is called {\it
non-decreasing} if $w_{(j-1)}\le w_{(j)}$ for all $j=1,\dotsc,
n$, so $J^-_\w=\emptyset$. It is called {\it reduced} if
$w_{(j-1)}< w_{(j)}$ for all $j=1,\dotsc, n$. Clearly, any
non-decreasing expression $\w$ for $w$ gives rise to a reduced
expression $\widehat{\w}$ of $w$ by discarding all $w_{(j)}$ with
$j\in J^\circ_\w$.
\end{defn}

The following definition is taken from {\cite[Definition
2.3]{Deo:Decomp}}.
\begin{defn}[Distinguished subexpressions]
Let $\w$ be a reduced expression for $w\in W$ with factors
$(s_{i_1},\dotsc, s_{i_n})$. Let $v\le w$. Then by a {\it
subexpression} for $v$ in $\w$, we mean an expression
$\v=(v_{(0)},v_{(1)},v_{(2)},\dotsc, v_{(n)})$ such that
\begin{equation*}
v_{(j)}\in\left\{v_{(j-1)},\, v_{(j-1)}s_{i_j}\right\} \qquad
\text{for all $j=1,\dotsc n$,}
\end{equation*}
and $v_{(n)}=v$. In particular there is always the ``empty''
subexpression $(1,\dotsc, 1)$ for $1$.

A subexpression $\v$ of $\w$ as above is called {\it
distinguished} if we have
\begin{equation}\label{e:dist}
v_{(j)}\le v_{(j-1)}\ s_{i_j}\quad \text{for all
$j\in\{1,\dotsc,n\}$}.
\end{equation}
In other words, if right  multiplication by $s_{i_j}$ decreases
the length of $v_{(j-1)}$, then  in a distinguished subexpression
the component $v_{(j)}$ must be given by
$v_{(j)}=v_{(j-1)}s_{i_j}$.

We write $\v\prec\w$ if $\v$ is a distinguished subexpression of
$\w$.
\end{defn}

\subsection*{Examples} For $w=w_0$ in $A_3$ and the reduced
expression $\w$ with factors $(s_3, s_2, s_1, s_3, s_2, s_3)$,
the only distinguished subexpression for $s_2 s_3$ is
\begin{align}\label{e:pos1}
\v=(1,1,1,1,1,s_2,s_2 s_3).
\end{align}
In particular, the subexpression indicated in \eqref{e:subexpr} is
not distinguished. If $v=s_2$, then we have four distinguished
subexpressions for $v$ in $\w$,
\begin{align}\label{e:pos2}
 \v&=(1,1,1,1,1,s_2,s_2), &\quad
 (s_3\, s_2 \, s_1\, s_3
 \, \underline{s_2} \,s_3), \\
 \v&=(1,s_3,s_3,s_3,1,s_2,s_2), &\quad
 (\underline{s_3}\, s_2 \, s_1\,\underline {s_3} \, \underline{s_2}
 \,s_3),\\
 \v&=(1,s_2,s_2,s_2 s_3,s_2 s_3,s_2), &
 \quad (s_3\, \underline{s_2} \, s_1\,\underline{s_3}
 \,s_2 \,\underline{s_3}),  \\
\v&=(1,s_3,s_3 s_2,s_3 s_2,s_3 s_2 s_3,s_2 s_3,s_2), & \quad
(\underline{s_3}\, \underline{s_2} \, s_1\,\underline{s_3}
 \,\underline{s_2} \,\underline{s_3}).
\end{align}

\begin{defn}[Positive subexpressions]
Let $\w$ be a reduced expression with factors $(s_{i_1},\dotsc,
s_{i_n})$. We call a subexpression $\v$ of $\w$ {\it positive}
if
 \begin{equation}\label{e:PositiveSubexpression}
v_{(j-1)}< v_{(j-1)}s_{i_j}
 \end{equation}
for all $j=1,\dotsc,n$.
\end{defn}

Note that \eqref{e:PositiveSubexpression} is equivalent to
$v_{(j-1)}\le v_{(j)}\le v_{(j-1)}s_{i_j}$. So in other words a
positive subexpression is one that is distinguished and
non-decreasing. In the examples above only \eqref{e:pos1} and
\eqref{e:pos2} are positive.
\begin{lem}\label{l:positive}
Given $v\le w$ in $W$ and a reduced expression $\w$ for $w$, then
there is a unique positive subexpression $\v_+$ for $v$ in $\w$.
\end{lem}
\begin{proof}
We construct $\v_+=(v_{(0)},\dotsc,v_{(n)})$ starting from the
right with $v_{(n)}=v$. The inequality $v_{(j-1)}<
v_{(j-1)}s_{i_j}$ says that $v_{(j-1)}$ cannot have a reduced
expression ending in $s_{i_j}$. If $v_{(j)}$ has such a reduced
expression then we must set $v_{(j-1)}=v_{(j)}s_{i_j}$. If
$v_{(j)}$ does not, then $v_{(j-1)}=v_{(j)}$. To summarize,
$v_{(j-1)}$ is given by
\begin{equation*}
v_{(j-1)}=\begin{cases} v_{(j)}s_{i_j}&\text{ if
$v_{(j)}s_{i_j}<v_{(j)}$,}\\
v_{(j)} &\text{ otherwise.}
\end{cases}
\end{equation*}
This along with $v_{(n)}=v$ clearly defines (uniquely) the desired
positive subexpression of $\w$.
\end{proof}
The positive subexpression $\v_+$ is in a sense the right-most
subexpression for $v$ in $\w$ that is non-decreasing.

\section{Deodhar's decomposition}
\label{s:parameterizations}

\subsection{Bruhat decomposition}\label{s:bruhatdecomposition}
Let us identify the flag variety
with the variety $\B$ of Borel subgroups, via
 $$
 gB^+\ \longleftrightarrow\ g\cdot B^+:=g B^+ g\inv.
 $$
We have the Bruhat decompositions,
 $$
 \B=\bigsqcup_{w\in W}B^+\dot w\cdot B^+=\bigsqcup_{w\in W}B^-\dot w\cdot
 B^-,
 $$
of $\B$ into $B^+$-orbits called {\it Bruhat cells}, and
$B^-$-orbits called {\it opposite Bruhat cells}. Let
$\alpha^1,\dotsc, \alpha^n$ be the positive roots made negative by
$w\inv$. Recall that the Bruhat cell $B^+\dot w\cdot B^+$ can be
identified with the product of root subgroups
\begin{equation}\label{e:rootsubgroups}
U^+\cap \dot w U^-\dot w\inv=U^+_{\alpha^1} U^+_{\alpha^2}\cdots
U^+_{\alpha^n}\cong \Ck^n
\end{equation}
via $u\mapsto u\dot w\cdot B^+$. Moreover,
\begin{equation*}
 U^+=(U^+\cap \dot w U^-\dot w\inv)(U^+\cap \dot w U^+\dot w\inv),
\end{equation*}
where the second factor is a product of the remaining positive
root subgroups. Given a reduced expression
$\w=(w_{(0)},w_{(1)},\dotsc, w_{(n)})$ with factors
$(s_{i_1},\dotsc, s_{i_n})$, the positive roots sent to negative
roots by $w\inv$ can be listed as
\begin{equation}\label{e:PosToNeg}
\alpha^1_{\w}=\alpha_{i_1}\ ,\quad \alpha^2_\w= w_{(1)}\cdot
\alpha_{i_2}\ ,\quad \alpha^3_\w=w_{(2)}\cdot\alpha_{i_3}\ ,\
\dotsc\dotsc\ ,\quad \alpha^n_\w= w_{(n-1)}\cdot \alpha_{i_n}.
\end{equation}
Therefore another way to write down the parametrization of the
Bruhat cell $B^+\dot w\cdot B^+$ is by
\begin{equation} \k^{n}\overset\sim\To B^+\dot w\cdot B^+~:\quad
(m_1,\dotsc, m_n)\mapsto x_{i_1}(m_1)\dot s_{i_1}\cdots
x_{i_n}(m_n)\dot s_{i_n}\cdot B^+.
\end{equation}
If one moves all the simple reflections to the right (conjugating
the intermediate simple root subgroups), then what remains on the
left is a product of root subgroups corresponding to precisely the
roots listed in \eqref{e:PosToNeg}.

\subsection{Relative position}\label{s:RelativePosition}
Consider the product $\B\x \B$ with $G$ acting diagonally. Let
$B_1=g_1\cdot B^+$ and $B_2=g_2\cdot B^+$. Then there is a unique
$w\in W$ such that ${g_1}\inv g_2\cdot B^+\in B^+\dot w\cdot
B^+$. Equivalently, $w$ is the unique Weyl group element such that
 $$
 (B_1,B_2)\, \in \text{ $G$-orbit of } (B^+,\dot w\cdot B^+).
 $$
We call $w$ the {\it relative position} of $(B_1,B_2)$ and write
 $$B_1\overset w\To B_2.$$
For example $B_1\overset 1\To B_2$ implies $B_1=B_2$. And
$B^+\overset {w}\To B$ says that $B$ lies in the Bruhat cell
$B^+\dot w\cdot B^+$. While $B^-\overset {w}\To B$ means that $B$
lies in the opposite Bruhat cell $B^-\dot w_0\dot w\dot w_0 \cdot
B^-$. We will also use the notation
 $$
 (B_1, B_2)\sim (B_1',B_2')
 $$
to indicate that $(B_1, B_2)$ and $(B_1', B_2')$ in $\B\x\B$ are
conjugate under $G$.

The following assertions follow from the definitions and standard
properties of the Bruhat decomposition.
\begin{enumerate}
 \item If $B_1\overset w\To B_2$ and $g\in G$, then also $g\cdot
 B_1\overset w\To g\cdot B_2$.
 \item If $B_1\overset s\To B_2\overset s\To B_3$ for a simple reflection $s$, then
$B_1\overset s \To B_3$ or $B_1=B_3$.
 \item If $ B_1\overset v\To B_2\overset w\To B_3$ and
$\ell(vw)=\ell(v)+\ell(w)$, then $B_1\overset{vw}\To B_3$.
 \item If $B_1\overset w\To B_2$, then $B_2\overset{w\inv}\To B_1$.
\end{enumerate}
We will make use of these properties freely.

\subsection{Reduction maps} \label{s:reductionmaps}

Suppose $w=vv'$ with $\ell(w)=\ell(v)+\ell(v')$. Then the set of positive
roots sent to negative roots by $v\inv$ is a subset of the
positive roots made negative by $w\inv$, by \eqref{e:PosToNeg}.
Under these circumstances one can define a morphism
\begin{align*}
 \pi^w_{v}: B^+\dot w\cdot B^+& \to B^+\dot v\cdot B^+\\
   \quad   b\dot w\cdot B^+ &\mapsto b\dot v\cdot B^+,
\end{align*}
where $b\in B^+$. The map $\pi^w_v$ is well-defined since $B^+\cap
\dot w B^+\dot w\inv\ \subseteq\  B^+\cap \dot v B^+ \dot v\inv$.
For a given $B\in B^+\dot w\cdot B^+$, the element $\pi^w_{v}(B)$
is characterized by the property
 \begin{equation}\label{e:reduction}
 B^+\ \overset{v}\To\ \pi^w_{v}(B)\ \overset {v\inv w}\To \ B.
 \end{equation}
Let us call $\pi^w_{v}$ a reduction map.

\subsection{Deodhar's theorem}\label{s:DeoComps}
\begin{defn}\label{d:intermediate}
For $v,w\in W$ define
\begin{equation*}
\mathcal R_{v,w}:=B^+\dot w\cdot B^+\cap B^-\dot v\cdot
B^+=\{B\in\B\ |\ B^+\overset {w}\To B\overset{w_0
v}\longleftarrow B^- \}.
\end{equation*}
\end{defn}
The intersection $\mathcal R_{v,w}$ is non-empty precisely if
$v\le w$. And in that case Kazhdan and Lusztig proved that over an
algebraically closed field it is
irreducible of dimension $\ell(w)-\ell(v)$, see
\cite{KaLus:Hecke}~\S 1. If $v=w$ then $\mathcal R_{w,w}=\{\dot
w\cdot B^+\}$.

Suppose now that $\w$ is a reduced expression for $w\in W$ with factors
$(s_{i_1},\dotsc, s_{i_n})$, and $B\in\mathcal R_{v,w}$. Using the
reduction maps we can associate to $B$ uniquely a sequence of
`intermediate' Borel subgroups
 $$
B^+=B_0\overset{s_{i_1}}\To B_1\overset{s_{i_2}} \To
B_2\overset{s_{i_3}}\To\cdots\overset{s_{i_n}}\To B_n=B,
 $$
where $B_{k}=\pi^{w}_{w_{(k)}}(B)$. By construction
$B^+\overset{w_{(k)}}\To B_k$. However, the position of $B_k$ with
respect to $B^-$, or the opposite Bruhat cell containing $B_k$,
defines a new element $v_{(k)}\in W$ by
 $$
 B_k\in B^- v_{(k)}\cdot B^+.
 $$
For $\w$ as above and a sequence $\v:=(v_{(0)},\dotsc, v_{(n)})$
we define the {\it Deodhar component} $\mathcal R_{\v,\w}$ in $
\mathcal B$ by
\begin{equation}\label{e:DeoDef}
\begin{aligned}
\mathcal R_{\v,\w}&:=\{B\in\mathcal R_{v,w}\ |\
\pi^{w}_{w_{(k)}}(B)\in B^- v_{(k)}\cdot B^+\ \}\\
 &\ =\{B\in\mathcal R_{v,w}\ |\
\pi^{w}_{w_{(k)}}(B)\in \mathcal R_{v_{(k)},w_{(k)}}\ \}.
\end{aligned}
\end{equation}

\begin{thm}[\cite{Deo:Decomp}~Theorem~1.1]\label{t:deodhardecomposition}
Suppose $w\in W$ and $B\in B^+\dot w\cdot B^+$, and fix a reduced
expression $\w=(w_{(0)},w_{(1)},\dotsc, w_{(n)})$ for $w$.
\begin{enumerate}
\item The Deodhar component $\mathcal R_{\v,\w}$ is
nonempty if and only if $\v$
is a distinguished subexpression of $\w$.
\item  If $\v\prec\w$, then $\mathcal R_{\v,\w}\cong
(\Ck^*)^{|J^{\circ}_\v|}\x \Ck^{|J^-_\v|}$, where $J^{\circ}_\v$
and $J^-_\v$ are as in Definition~\ref{d:Js}.
\end{enumerate}
\end{thm}
Another proof of this theorem will be contained in the next
section. If the reduced expression $\w$ is fixed, then as a
corollary of the theorem one has a decomposition
 \begin{equation}\label{e:Rdecomp}
 \mathcal R_{v,w}=\bigsqcup_{\v} \mathcal R_{\v,\w},
 \end{equation}
where the union is over all distinguished subexpressions for $v$
in $\w$. Note that the Deodhar component $\mathcal R_{\v_+,\w}$
corresponding to the
unique positive subexpression for $v$ in $\w$ has dimension
$|J^\circ_{\v_+}|=\ell(w)-\ell(v)$. So if $\k$ is algebraically closed
then it is dense in
$\mathcal R_{v,w}$. This also holds for $\k=\R$ since
$\mathcal R_{\v_+,\w}(\R)$ is Zariski dense in $\mathcal R_{\v_+,\w}(\C)$.
Finally, for $\k=\R$ or $\k=\C$, it holds that
$\mathcal R_{\v_+,\w}$ is open dense
in $\mathcal R_{v,w}$
with respect to the usual Hausdorff topology.

Suppose we fix a reduced expression $\w_0$ for the longest
element $w_0$. Then for any $w\in W$ the positive subexpression
for $w$ in $\w_0$ determines a reduced expression $\widehat\w_+$
for $w$. Therefore we have a decomposition of the whole flag
variety,
\begin{equation}\label{e:DeoDecomp}
 \mathcal B=\bigsqcup_{w\in W}\left(\bigsqcup_{\v\prec \widehat\w_+} \mathcal
 R_{\v,\widehat\w_+}\right),
\end{equation}
which we may call the Deodhar decomposition of $\B$ corresponding
to $\w_0$.

\begin{rem}
The varieties $\mathcal R_{v,w}$ may be defined over a finite
field $\k=\mathbb F_q$. In this setting the number of points
determine the $R$-polynomials $R_{v,w}(q)=\#(\mathcal
R_{v,w}(\mathbb F_q))$ introduced by Kazhdan and Lusztig
\cite{KaLus:Hecke} to give a recursive formula for the
Kazhdan-Lusztig polynomials. This is the origin of the notation
$\mathcal R_{v,w}$ as well as Deodhar's original application of
the theorem.
The decompositions \eqref{e:Rdecomp} together with the
isomorphisms $\mathcal R_{\v,\w}(\mathbb F_q)\cong (\mathbb
F_q^*)^{|J_\v^{\circ}|}\x{\mathbb F_q}^{|J_{\v}^-|}$ give formulas
for the $R$-polynomials.
\end{rem}

\section{Explicit parametrizations of Deodhar components}
\label{s:param} Let $\w$ be a reduced
expression with factors $(s_{i_1},\dotsc, s_{i_n})$, and
$\v\prec\w$.
\begin{defn} \label{d:factorization}
Define a subset $G_{\v,\w}$ in $G$ by
\begin{equation}\label{e:Gvw}
G_{\v,\w}=\left\{g= g_1 g_2\cdots g_n \left
|\begin{array}{ll}
 g_k= x_{i_k}(m_k)\dot s_{i_k}\inv& \text{ if $k\in J^-_\v$,}\\
 g_k= y_{i_k}(t_k)& \text{ if $k\in J^{\circ}_\v$,}\\
 g_k=\dot s_{i_k}& \text{ if $k\in J^+_\v$,}
 \end{array}\quad \text{
for $t_k\in\Ck^*,\, m_k\in\Ck$. }\right. \right\}.
\end{equation}
There is an obvious map $(\Ck^*)^{J^{\circ}_\v}\x\Ck^{J^-_\v}\to
G_{\v,\w}$ defined by the parameters $t_k$ and $m_k$ in
\eqref{e:Gvw}. For $\v =\w=(1)$ we define $G_{\v,\w}=\{1\}$.

\end{defn}
The following proposition gives an explicit parametrization for
the Deodhar component $\mathcal R_{\v,\w}$.
\begin{prop}\label{p:parameterization}
The map $(\Ck^*)^{J^{\circ}_\v}\x\Ck^{J^-_\v}\to G_{\v,\w}$ from
Definition~\ref{d:factorization} is an isomorphism. The set
$G_{\v,\w}$ lies in $U^-\dot v\cap B^+\dot w B^+$, and the
assignment $g\mapsto g\cdot B^+$ defines an isomorphism
\begin{align}\label{e:parameterization}
G_{\v,\w}&\overset\sim\To \mathcal R_{\v,\w}.
\end{align}
\end{prop}
The special case $v=1$ and $w=w_0$ of this proposition already
appears in \cite{MarRie:G2}~Proposition~2.5. The proof below is
analogous to the one we gave for that special case, and of course
also similar to Deodhar's proof of
Theorem~\ref{t:deodhardecomposition}.(2), although his is
ultimately a different isomorphism.
\begin{proof}
Let $\w=(w_{(0)},\dotsc, w_{(n)})$ be a reduced expression with
factors $(s_{i_1},\dotsc, s_{i_n})$, and $\v=(v_{(0)},\dotsc,
v_{(n)})$. The proof is by induction on $n$. If $n=0$ then
$\v=\w=(1)$ and the isomorphism \eqref{e:parameterization} is the
trivial one $1\mapsto B^+$. There is nothing more to check. For
$n>0$ let $\w':=(w_{(0)},\dotsc, w_{(n-1)})$ and similarly
$\v'=(v_{(0)},\dotsc, v_{(n-1)})$, the truncations of $\v$ and $\w$.
Also set $w'=w_{(n-1)}$ and $v'=v_{(n-1)}$. We may assume the
proposition is true for $\v',\w'$.

It is easy to check that $G_{\v',\w'}\x \mathbb \Ck
\overset\sim\To (\pi_{w'}^{w})\inv(\mathcal R_{\v',\w'})$ via the
map $(g',m)\mapsto g'x_{i_{n}}(m)\dot s_{i_{n}}\cdot B^+$, using
for example \eqref{e:reduction} and properties of relative
position. And we have a commutative diagram
 \begin{equation}\label{e:CD}
 \begin{array}{ccc}
 G_{\v',\w'}\x \mathbb \Ck &\overset\sim\To &
  (\pi_{w'}^{w})\inv(\mathcal R_{\v',\w'})\\
   pr_1\downarrow & & \downarrow \pi_{w'}^{w}\\
 G_{\v',\w'}&\overset\sim\To &\mathcal R_{\v',\w'}.
  \end{array}
 \end{equation}
Now let $B\in (\pi_{w'}^{w})\inv(\mathcal R_{\v',\w'})$, so
$B=g'x_{i_{n}}(m)\dot s_{i_{n}}\cdot B^+$ for some $g'\in
G_{\v',\w'}$ and $m\in\Ck$.
We consider two cases.
\begin{enumerate}
\item[$(i)$] Suppose $m=0$. Then
$B=g'\dot s_{i_n}\cdot B^+$. Since $g'\in U^-\dot v'$ we have
$B\in B^-\dot v'\dot s_{i_n}\cdot B^+$.
\item[$(ii)$] Suppose $m\ne 0$.
Then the identity \eqref{e:siIdentity} implies
$x_{i_n}(m)\dot s_{i_n}\cdot B^+=y_{i_n}(m\inv)\cdot B^+$. So we
may write $B$ in two different ways,
 $$B=g'x_{i_n}(m)\dot s_i\cdot B^+=g'y_{i_n}(m\inv)\cdot B^+.$$
\begin{itemize}
\item
If $v's_{i_n}>v'$, then $\dot v' y_{i_n}(m\inv){\dot v'}{}\inv\in
U^-$. In this case we have
 $$
g:=g'y_{i_n}(m\inv)\in U^-\dot v'
 $$
and $B=g\cdot B^+\in B^-\dot v'\cdot B^+$.
\item
If $v' s_{i_n}< v'$, then $\dot v'x_{i_n}(m){\dot v'}{}\inv\in
U^-$. Therefore we have
 $$
 g:=g'x_{i_n}(m)\dot s_{i_n}\inv\in U^-\dot v'\dot
s_{i_n}\inv
 $$
and $B=g\cdot B^+\in B^-\dot v'\dot s_{i_n}\inv\cdot B^+$.
\end{itemize}
\end{enumerate}

Note that in both cases, $(i)$ and $(ii)$, if $v's_{i_n}<v'$ we
have $B\in B^-\dot v'\dot s_{i_n}\cdot B^+$. This explains
Theorem~\ref{t:deodhardecomposition}.(1). We now use the above to
analyze the possibilities for an element $B\in\mathcal
R_{\v,\w}\subseteq(\pi^{w'}_w)\inv(\mathcal R_{\v',\w'})$ and
complete the proof of the proposition.
\begin{enumerate}
\item[{$(1)$}]
 Suppose $n\in J^-_\v$. Then both $(i)$ and $(ii)$
are possible. Therefore $\mathcal
R_{\v,\w}=(\pi_{w}^{w'})\inv(\mathcal R_{\v',\w'})$, and we have
 $$
G_{\v,\w}\cong G_{\v',\w'}\x \k\overset \sim\To \mathcal R_{\v,\w},
 $$
via $g_1\dotsc g_{n-1}x_{i_n}(m_n)\dot s_{i_n}\inv\mapsto
(g_1\dotsc g_{n-1},m_n)$ and \eqref{e:CD}.
\item[{$(2)$}]
If $n\in J^\circ_\v$ then $v_{(n)}=v_{(n-1)}$ so only case $(ii)$
is possible. Then \eqref{e:CD} restricts to give
 $$
 G_{\v,\w}\cong G_{\v',\w'}\x\k^*\overset\sim\To\mathcal R_{\v,\w},
 $$
where the identification $G_{\v,\w}\cong G_{\v',\w'}\x \k^*$ is
via $g_{1}\dotsc g_{n-1}y_{i_n}(t_n)\mapsto (g_{1}\dotsc
g_{n-1},t_n\inv)$.
\item[$(3)$]
Finally, if $n\in J^+_\v$ then only case $(i)$ is possible and
$B=g'\dot s_{i_n}\cdot B^+$. Therefore \eqref{e:CD} induces
 $$
G_{\v,\w}\cong G_{\v',\w'}\x\{\, 1\,\}\overset\sim\To R_{\v,\w}.
 $$
\end{enumerate}
In each case $G_{\v,\w}\subset U^-\dot v_{(n)}$, where we note
that in $(1)$ above, $v_{(n)}=v's_{i_n}<v'$ implies $\dot
v_{(n)}=\dot v'\dot s_{i_n}\inv$. The inclusion $G_{\v,\w}\subset
B^+\dot w B^+$ is clear.
\end{proof}

\begin{rem} Let $\w_0$ be a fixed reduced expression for $w_0$. Then
Deodhar's decomposition \eqref{e:DeoDecomp} of $\mathcal B$ together with
Proposition~\ref{p:parameterization} gives rise to a set theoretic crossection
 $$
\phi:\mathcal B\To G,
 $$
defined on each Deodhar component
$\mathcal R_{\v,\widehat\w_+}\subset \mathcal B$
as the inverse
of $G_{\v,\widehat\w_+}\overset\sim\To \mathcal R_{\v,\widehat\w_+}$.
To describe the map
$\phi$ more explicitly
we must, firstly, explain how to determine the Deodhar component of an
element of $\mathcal B$ and, secondly, give formulas for the individual maps
$\mathcal R_{\v,\w}\to G_{\v,\w}$.
\end{rem}

\section{Deodhar components in terms of minors}\label{s:algorithm}
Suppose $B$ lies in a particular Bruhat cell,
$B=z\dot w\cdot B^+$ for $z\in U^+$.
In this section we determine the conditions on $z$ for $B$
to lie in a Deodhar component $\mathcal
R_{\v,\w}$. The conditions will be expressed in terms of (generalized)
minors of $z$.

Let
$V(\lambda)$ be the Weyl module of $G$ with highest weight
$\lambda$. In the following $\lambda$ will often be a fundamental
weight $\omega_i$. Consider the weight space decomposition
$V(\lambda)=\oplus_{\mu} V(\lambda)_{\mu}$, and denote by
$pr_{\mu}:V(\lambda)\to V(\lambda)_{\mu}$ the corresponding
projections. Let us fix a highest weight vector $\xi_{\lambda}$.
Then the element $\dot w\cdot \xi_\lambda\in V(\lambda)$ for $w\in
W$ spans the extremal weight space $V(\lambda)_{w\lambda}$. In
this way, the choice of highest weight vector gives rise to a
canonical choice of basis vectors for all the extremal weight
spaces.
\begin{lem}
If $w\lambda=w'\lambda $, then $\dot
w\cdot\xi_\lambda=\dot w'\cdot\xi_\lambda$.
\end{lem}
\begin{proof}
It is necessary only to check that $\dot v\cdot\xi_\lambda=\xi_\lambda$
whenever $v\lambda=\lambda$. Since the stabilizer of $\lambda$ is a parabolic
subgroup of $W$ we may assume $v$ is a simple reflection $s_i$. Then
$\dot s_i=x_i(-1)y_i(1)x_i(-1)$ and the statement is clear.
\end{proof}

\begin{defn}[Generalized minors] \label{d:minors}
For $\eta \in V(\lambda)$ define $\left< \eta,\dot
w\cdot\xi_{\lambda} \right>$ to be the coefficient in $\eta$ of the extremal
weight vector $\dot w\cdot \xi_\lambda$. That is, with notation as
above,
 $$
pr_{w\lambda}(\eta)=\left< \eta,\dot w \cdot\xi_{\lambda} \right>\
\dot w\cdot \xi_{\lambda}.
 $$
For two extremal weights $w\lambda$ and $w'\lambda$ we then have a
regular function $\Delta_{w'\lambda}^{w\lambda}$ on $G$ defined by
 $$
 \Delta_{w'\lambda}^{w\lambda}(g):=\left< g\dot w'\cdot
\xi_{\lambda}\, ,\,\dot w \cdot\xi_{\lambda}\right>.
 $$
Since any weight lies in the Weyl group orbit of a unique
dominant weight, this notation is unambiguous.

It is not hard to see that
$\Delta_{w'\lambda}^{w\lambda}$ coincides with the regular function
$\Delta_{w\lambda, w'\lambda}$ defined in~\cite[Definition 1.4]
{FoZe:DoubleCells}.

The functions  $\Delta^{w\omega_i}_{w'\omega_i}$, where $\omega_i$
ranges through the set of fundamental weights, are called {\it
minors} or {\it generalized minors}. If $G=SL_d$ with the standard
pinning then $\Delta^{w\omega_i}_{w'\omega_i}$ is precisely the
usual $i\x i$ minor, where $w\omega_i$ encodes the row set and
${w'\omega_i}$ the column set.
\end{defn}

\begin{defn}[Chamber minors]  \label{d:chamberminors}
Suppose $\w$ is a reduced expression and $\v\prec \w$ a
distinguished subexpression.

\begin{enumerate}\item
The minors $\Delta_{w_{(k)}\omega_{i_k}}^{v_{(k)}\omega_{i_k}}$
for $k=0,1,\ldots ,n$ are called the {\em standard chamber minors}
for $\v$ and $\w$.
\item
The minors $\Delta_{w_{(k)}\omega_{i_k}}^{v_{(k-1)}\omega_{i_k}}$
for $k\in J^-_\v\cup J^+_\v$ are called the {\em special chamber
minors} for $\v$ and $\w$.
\end{enumerate}
Note that $\Delta_{w_{(k)}\omega_{i_k}}^{v_{(k-1)}\omega_{i_k}}=
\Delta_{w_{(k)}\omega_{i_k}}^{v_{(k)}\omega_{i_k}}$ if $k\in
J^\circ_\v$.
\end{defn}

\begin{prop}\label{p:algorithm} Let $B=z\dot w\cdot B^+$ for
$z\in U^+$, and $\mathbf{w}$ be a reduced expression with factors
$(s_{i_1},s_{i_2},\ldots ,s_{i_n})$. Then $B$ lies in the the
Deodhar component $\mathcal R_{\v,\w}$ where
$\v=(v_{(0)},v_{(1)},\dotsc, v_{(n)})$ is constructed as follows.
Let $v_{(0)}=1$. Suppose that $k\geq 1$ and $v_{(k-1)}$ has
already been defined. \begin{itemize}
\item[$(a)$] If $v_{(k-1)}s_{i_k}>v_{(k-1)}$ and
$\Delta_{w_{(k)}\omega_{i_k}}^{v_{(k-1)}\omega_{i_k}}(z)\not=0$,
then $v_{(k)}=v_{(k-1)}$.
\item[$(b)$] If $v_{(k-1)}s_{i_k}>v_{(k-1)}$ and
$\Delta_{w_{(k)}\omega_{i_k}}^{v_{(k-1)}\omega_{i_k}}(z)=0$, then
$v_{(k)}=v_{(k-1)}s_{i_k}$.
\item[$(c)$] If $v_{(k-1)}s_{i_k}<v_{(k-1)}$, then
$v_{(k)}=v_{(k-1)}s_{i_k}$.
\end{itemize}
\end{prop}
\begin{rem}
Note that in the situation of the proposition
$\Delta_{w_{(k)}\omega_{i_k}}^{v_{(k)}\omega_{i_k}}(z)\ne 0$ for
all $k=1,\dotsc, n$, as follows from the definition of the
$v_{(k)}$. The chamber minors give rise to well-defined maps
(which we denote in the same way),
\begin{equation*}
\begin{array}{rclll}
\Delta_{w_{(k)}\omega_{i_k}}^{v_{(k)}\omega_{i_k}}:&
\mathcal R_{\v,\w}\to \k^*: & z\dot w\cdot B^+\mapsto
\Delta^{v_{(k)}\omega_{i_k}}_{w_{(k)}\omega_{i_k}}(z)&\qquad k=1,\dots, n,\\
\Delta^{v_{(k-1)}\omega_{i_k}}_{w_{(k)}\omega_{i_k}}:&
\mathcal R_{\v,\w}\to \k~: & z\dot w\cdot B^+\mapsto
\Delta^{v_{(k-1)}\omega_{i_k}}_{w_{(k)}\omega_{i_k}}(z)&\qquad k\in J^{-}_\v.
\end{array}
\end{equation*}
\end{rem}
\begin{proof}
By Theorem~\ref{t:deodhardecomposition}.(1) we have that $\v$ is a
distinguished subexpression of $\w$. Therefore $(c)$ holds. Now
suppose $v_{(k-1)}s_{i_k}>v_{(k-1)}$. We have either
\begin{enumerate}
\item
$v_{(k)}=v_{(k-1)}$ and $z\dot w_{(k)}\cdot B^+\in B^-\dot
v_{(k-1)}\cdot B^+$, or
\item
$v_{(k)}=v_{(k-1)}s_{i_k}$ and $z\dot w_{(k)}\cdot B^+\in
B^-\dot v_{(k-1)}\dot s_{i_k}\cdot B^+$.
\end{enumerate}
We can distinguish between these two cases by looking just at the
representation $V_{\omega_{i_k}}$. In the first case, the highest
weight occurring in $z\dot w_{(k)}\cdot \xi_{\omega_{i_k}}$ is
$v_{(k-1)}\omega_{i_k}$, and hence
$\Delta_{w_{(k)}\omega_{i_k}}^{v_{(k-1)}\omega_{i_k}}(z)\not=0$.
In the second case, the highest weight occurring in $z\dot
w_{(k)}\cdot \xi_{\omega_{i_k}}$ is
$v_{(k-1)}s_{i_k}\omega_{i_k}$, which is lower than
$v_{(k-1)}\omega_{i_k}$ since $v_{(k-1)}s_{i_k}>v_{(k-1)}$.
Therefore we have
$\Delta_{w_{(k)}\omega_{i_k}}^{v_{(k-1)}\omega_{i_k}}(z)=0$.
\end{proof}

As a reformulation of Proposition~\ref{p:algorithm} we have the
following description of $\mathcal{R}_{\mathbf{v},\mathbf{w}}$
inside the Bruhat cell $B^+\dot w\cdot B^+$.

\begin{cor}\label{c:DeoInequalities}
Suppose $\w$
is a reduced expression of $w$ and $\v\prec \w$ a distinguished
subexpression, with $J^+_\v$ and $J^{\circ}_\v$ as in
Definition~\ref{d:Js}. Then $\mathcal{R}_{\mathbf{v},\mathbf{w}}$
may be described by
 $$
\mathcal{R}_{\mathbf{v},\mathbf{w}}=\left\{z\dot w\cdot B^+\,
\left |\, z\in U^+\, ;\,
\begin{array}{ll}\Delta_{w_{(k)}\omega_{i_k}}^{v_{(k-1)}\omega_{i_k}}(z)=
 0 & \text{for all $k\in J^+_\v$,}\\
 \Delta_{w_{(k)}\omega_{i_k}}^{v_{(k)}\omega_{i_k}}(z)\ne 0
 & \text{for all $k\in J^{\circ}_\v$}
 \end{array}\right . \right\}.
 $$
\end{cor}
\section{The generalized Chamber Ansatz}
By Proposition~\ref{p:parameterization} a Deodhar component $\mathcal
R_{\v,\w}$
comes with isomorphisms
   \begin{equation}\label{e:par}
   (\k^*)^{J^\circ_\v}\x \k^{J^-_\v}\overset \sim \to
G_{\v,\w}\overset\sim\to \mathcal R_{\v,\w}.
   \end{equation}
The aim of this section is to describe an inverse to \eqref{e:par}.
The following theorem generalizes the Chamber Ansatz of Berenstein and
Zelevinsky~\cite{BeZel:TotPos}.
\begin{thm} (Generalized Chamber Ansatz) \\ \label{t:ansatz}
Let $B=z\dot w\cdot B^+\in {\mathcal R}_{v,w}$, where $z\in U^+$,
$v,w\in W$ and $v\leq w$. Let $\mathbf{w}=(w_{(0)},w_{(1)},\ldots ,w_{(n)})$
be a reduced expression for $w$ with factors
$(s_{i_1},s_{i_2},\ldots ,s_{i_n})$. Then $B$ lies in a Deodhar component
${\mathcal R}_{\mathbf{v},\mathbf{w}}$, where $\mathbf{v}=
(v_{(0)},v_{(1)},\ldots ,v_{(n)})$ is a distinguished
subexpression for $v$ in $\mathbf{w}$.
By Proposition~\ref{p:parameterization}, there is
$g\in G_{\mathbf{v},\mathbf{w}}$ such that $B=g\cdot B^+$.
By Definition~\ref{d:factorization} we can write
$g=g_1g_2\cdots g_n\in U^-\dot v\cap B^-\dot w B^+$, where
$$g_k=\left\{\begin{array}{ll}
y_{i_k}(t_k) & k\in J^{\circ}_{\mathbf{v}}, \\
\dot s_{i_k} & k\in J^+_{\mathbf{v}}, \\
x_{i_k}(m_k) \dot s_{i_k}^{-1} & k\in J^-_{\mathbf{v}}.
\end{array}
\right.
$$
For each $k$, let $g_{(k)}=g_1g_2\cdots g_k$ denote the partial product.
Then the following hold.
\begin{enumerate}
\item For $k\in J^{\circ}_{\mathbf{v}}$, we have: $$t_k=\frac{
\prod_{j\not=i_k}\Delta_{w_{(k)}\omega_j}^{v_{(k)}\omega_j}(z)^{-a_{j,i_k}}
}{\Delta_{w_{(k)}\omega_{i_k}}^{v_{(k)}\omega_{i_k}}(z)\Delta_{w_{(k-1)}\omega_{i_k}}^{v_{(k-1)}\omega_{i_k}}(z)}$$
\item For $k\in J^-_{\mathbf{v}}$, we have:
$$m_k=\frac{
\Delta_{w_{(k)}\omega_{i_k}}^{v_{(k-1)}\omega_{i_k}}(z)\Delta_{w_{(k-1)}\omega_{i_k}}^{v_{(k-1)}\omega_{i_k}}(z)}
{\prod_{j\not=i_k}\Delta_{w_{(k)}\omega_j}^{v_{(k)}\omega_j}(z)^{-a_{j,i_k}}}
-\Delta^{v_{(k-1)}\omega_{i_k}}_{s_{i_k}\omega_{i_k}}(g_{(k-1)}).$$
\end{enumerate}
\end{thm}

\begin{rem}
It is easy to check that the minors appearing in Theorem~\ref{t:ansatz}(1)
are the standard chamber minors of Definition~\ref{d:chamberminors}(1).
The formula for the
$m_k$ also involves the special chamber minors, as well as a correction term,
$\Delta^{v_{(k-1)}\omega_{i_k}}_{s_{i_k}\omega_{i_k}}(g_{(k-1)})$,
which can be computed recursively. It is an open problem to find a
closed formula in terms of minors of $z$ for this correction term.
See Section~\ref{s:coordinates} for another interpretation of the formulas in
Theorem~\ref{t:ansatz}.
\end{rem}

In order to prove Theorem~\ref{t:ansatz}, we will rewrite the chamber
minors as minors of $g_{(k)}=g_1g_2\cdots g_k$, for $k=0,1,\ldots ,n$
(Lemma~\ref{l:crossminor}). We will then compute these minors
(Lemma~\ref{l:gminor}) and substitute these formulas back into
the expressions in Theorem~\ref{t:ansatz}(1) and (2), finally showing that
they reduce to the coefficients $t_k$ and $m_k$ as claimed.
\begin{lem} \label{l:partialproduct}
For $k\in\{0,1,\ldots ,n\}$ we have
$g_{(k)}\cdot B^+=z\dot w_{(k)}\cdot B^+$.
\end{lem}
\begin{proof}
Consider the sequence $B^+,g_{(1)}\cdot B^+,g_{(2)}\cdot B^+,
\ldots ,g_{(n)}\cdot B^+$. Then clearly
$$B^+\ \overset{s_{i_1}}\To\ g_{(1)}\cdot B^+\ \overset {s_{i_2}}\To \
g_{(2)}\cdot B^+\ \overset{s_{i_3}}\To\ \cdots \overset{s_{i_n}}\To\
g_{(n)}\cdot B^+.$$
It follows from Section~\ref{s:reductionmaps} that
$g_{(k)}\cdot B^+=\pi^w_{w_{(k)}}(B)=z\dot w_{(k)}\cdot B^+$ for
$k=0,1,\ldots ,n$ as required.
\end{proof}

We note that Lemma~\ref{l:partialproduct} gives two expressions for the
intermediate Borel subgroups $B_k=\pi^w_{w_{(k)}}(B)$. Suppose that $\lambda$
is a dominant weight. Let $L_k$ denote the line in the module $V(\lambda)$
stabilised by $B_k$.
Then, by the above, we have $L_k=<z\dot w_{(k)}\cdot \xi_{\lambda}>
=<g_{(k)}\cdot \xi_{\lambda}>$, where $<\eta>$ denotes the line
spanned by $\eta\in V(\lambda)$. We use this fact in the following lemma
to compute minors of $z$ in terms of minors of $g_{(k)}$.

\begin{lem} \label{l:crossminor}
Let $\lambda$ be a dominant weight. Then
\begin{enumerate}
\item For $k\in \{0,1,\ldots ,n\}$, we have
$$\Delta_{w_{(k)}\lambda}^{v_{(k)}\lambda}(z)=
\frac{1}{\Delta_{\lambda}^{w_{(k)}\lambda}(g_{(k)})},\mbox{\ \ and}$$
\item For $k\in \{1,\ldots ,n\}$, we have
$$\Delta_{w_{(k)}\lambda}^{v_{(k-1)}\lambda}(z)=
\frac{\Delta_{\lambda}^{v_{(k-1)}\lambda}(g_{(k)})}
{\Delta_{\lambda}^{w_{(k)}\lambda}(g_{(k)})}.$$
\end{enumerate}
\end{lem}
\begin{proof}
Let $L_k=<z\dot w_{(k)}\cdot \xi_{\lambda}>=<g_{(k)}\cdot \xi_{\lambda}>$
be the line in $V(\lambda)$ defined above.
Since $z\in U^+$, we have
$\langle z\dot w_{(k)} \cdot \xi_{\lambda}, \dot w_{(k)} \cdot
\xi_{\lambda}\rangle=1$.
Therefore
\begin{equation}\label{e:crossequation}
z\dot w_{(k)} \cdot \xi_{\lambda}=\left(\frac{1}{\left< g_{(k)} \cdot
\xi_{\lambda},
\dot w_{(k)} \cdot \xi_{\lambda}\right>}\right) g_{(k)}\cdot \xi_{\lambda}.
\end{equation}
Comparing coefficients of $\dot v_{(k-1)}\cdot \xi_{\lambda}$ on both sides,
(2) immediately follows, and
comparing coefficients of $\dot v_{(k)}\cdot \xi_{\lambda}$ on both sides,
we obtain
$$\Delta_{w_{(k)}\lambda}^{v_{(k)}\lambda}(z)=
\frac{\Delta_{\lambda}^{v_{(k)}\lambda}(g_{(k)})}
{\Delta_{\lambda}^{w_{(k)}\lambda}(g_{(k)})}.$$
However, by Proposition~\ref{p:parameterization}, $g_{(k)}\in U^-\dot v_{(k)}$,
so $\langle g_{(k)}\cdot \xi_{\lambda},\dot v_{(k)}\cdot \xi_{\lambda}\rangle=
\Delta_{\lambda}^{v_{(k)}\lambda}(g_{(k)})=1$, and (1) follows.
\end{proof}

We now compute the minors of $g_{(k)}$ from Lemma~\ref{l:crossminor}.

\begin{lem} \label{l:gminor}
Let $k\in \{0,1,\ldots ,n\}$ and let $\lambda$ be a dominant weight.
Then we have:
\begin{enumerate}
\item $$\Delta_{\lambda}^{w_{(k)}\lambda}(g_{(k)})=
\prod_{l=1,\ l\in J^{\circ}_{\mathbf{v}}}^k
t_l^{\left< s_{i_{l+1}}s_{i_{l+2}}\cdots s_{i_k}\lambda,\alpha_{i_l}^{\vee}
\right>}
\prod_{l=1,\ l\in J^-_{\mathbf{v}}}^k
(-1)^{\left< s_{i_{l+1}}s_{i_{l+2}}\cdots s_{i_k}\lambda,\alpha^{\vee}_{i_l}
\right>}.$$
\item If $k\in J^-_{\mathbf{v}}$, then
$$m_k=-\Delta_{\omega_{i_k}}^{v_{(k-1)}\omega_{i_k}}(g_{(k)})
-\Delta^{v_{(k-1)}\omega_{i_k}}_{s_{i_k}\omega_{i_k}}(g_{(k-1)}).$$
\end{enumerate}
\end{lem}
\begin{proof}
(1) We prove the result for $g=g_{(n)}=g_1g_2\cdots g_n$. The result
for arbitrary $k$ follows since $g_{(k)}$ is defined in terms of the
distinguished subexpression
$\mathbf{v}_{(k)}=(v_{(0)},v_{(1)},v_{(2)},\ldots ,v_{(k)})$ for
$v_{(k)}$ in $w_{(k)}$ in the same way that $g$ is defined in terms
of the reduced subexpression $\mathbf{v}$ for $v$ in $\mathbf{w}$.
For $l=1,2,\ldots ,n+1$, let $g^{(l)}=g_lg_{l+1}\cdots g_n$.
We note that $g^{(l)}\in B^+\dot w^{(l)}B^+$ for $l=1,2,\ldots ,n+1$,
where $w^{(l)}=\dot s_{i_l}\dot s_{i_{l+1}} \cdots \dot s_{i_n}$.
We prove, by reverse induction on $l$, that
$$\Delta_{\lambda}^{w^{(l)}\lambda}(g^{(l)})=
\prod_{j=l,\ j\in J^{\circ}_{\mathbf{v}}}^n
t_j^{\left< s_{i_{j+1}}s_{i_{j+2}}\cdots
s_{i_n}\lambda,\alpha^{\vee}_{i_j}\right>}
\prod_{j=l,\ j\in J^-_{\mathbf{v}}}^n
(-1)^{\left< s_{i_{j+1}}s_{i_{j+2}}\cdots
s_{i_n}\lambda,\alpha^{\vee}_{i_j}\right>}.$$
The start of the induction is clear. Suppose that the result holds
for $l+1$, i.e. for $g^{(l+1)}=g_{l+1}\cdots g_n$,
and consider $g^{(l)}=g_lg_{l+1}\cdots g_n$.
Since $g^{(l+1)}\in B^+\dot w^{(l+1)}B^+$,
$g^{(l+1)}\cdot \xi_{\lambda}$ is a linear combination of elements of
$V(\lambda)$ of weight $\mu\geq w^{(l+1)}\lambda$.

Case (I). Suppose that $l\in J^{\circ}_{\mathbf{v}}$, so that
$g_l=y_{i_l}(t_l)$.
Then, using that $w^{(l+1)}\lambda$ and $w^{(l)}\lambda$ are extremal weights,
and $w^{(l)}\lambda=s_{i_l}w^{(l+1)}\lambda \le w^{(l+1)}\lambda$
we have
\begin{eqnarray*}
\Delta_{\lambda}^{w^{(l)}\lambda}(g^{(l)}) & = &
\Delta_{\lambda}^{w^{(l)}\lambda}(y_{i_l}(t_l)g^{(l+1)}) \\
& = & \left< y_{i_l}(t_l)g^{(l+1)}\cdot \xi_{\lambda},\dot w^{(l)}
\cdot \xi_{\lambda} \right> \\
& = & \left< y_{i_l}(t_l) \dot w^{(l+1)}\cdot\xi_\lambda,
        \dot w^{(l)}\cdot\xi_{\lambda}
      \right>
      \left< g^{(l+1)}\cdot \xi_{\lambda}, \dot w^{(l+1)}\cdot
        \xi_{\lambda}
      \right>.
\end{eqnarray*}
By equation~\eqref{e:siIdentity},
\begin{eqnarray*}
\Delta_{\lambda}^{w^{(l)}\lambda}(g^{(l)})
& = &
\left< x_{i_l}(t_l\inv)\alpha_{i_l}^{\vee}(t_l\inv)\dot
s_{i_l}x_{i_l}(t_l\inv) \dot w^{(l+1)}\cdot \xi_{\lambda},
\dot w^{(l)}\cdot \xi_{\lambda} \right>
\left< g^{(l+1)}\cdot \xi_{\lambda},\dot w^{(l+1)}\xi_{\lambda}\right>
\\
& = &
t_l^{\left< s_{i_{l+1}}\cdots
s_{i_n}\lambda,\alpha^{\vee}_{i_l}\right>}
\Delta_{\lambda}^{w^{(l+1)}\lambda}(g^{(l+1)}).
\end{eqnarray*}

Case (II). Suppose that $l\in J^+_{\mathbf{v}}$, so that $g_l=\dot s_{i_l}$.
Then
\begin{eqnarray*}
\Delta_{\lambda}^{w^{(l)}\lambda}(g^{(l)}) & = &
\Delta_{\lambda}^{w^{(l)}\lambda}(\dot s_{i_l}g^{(l+1)}) \\
& = & \left< \dot s_{i_l}g^{(l+1)}\cdot \xi_{\lambda},\dot w^{(l)}
\cdot \xi_{\lambda} \right>.
\end{eqnarray*}
It is thus clear that
$$\Delta_{\lambda}^{w^{(l)}\lambda}(g^{(l)})=
\Delta_{\lambda}^{w^{(l+1)}\lambda}(g^{(l+1)}).$$

Case (III). Suppose that $l\in J^-_{\mathbf{v}}$, so that
$g_l=x_{i_l}(m_l)\dot s_{i_l}^{-1}$. Then, using the fact that
that $w^{(l+1)}\lambda\geq w^{(l)}\lambda=s_{i_l}w^{(l+1)}\lambda$
are extremal weights, we have:
\begin{eqnarray*}
\Delta_{\lambda}^{w^{(l)}\lambda}(g^{(l)}) & = &
\Delta_{\lambda}^{w^{(l)}\lambda}(x_{i_l}(m_l)\dot s_{i_l}^{-1}
g^{(l+1)}) \\
& = & \left< x_{i_l}(m_l)\dot s_{i_l}^{-1}g^{(l+1)}\cdot \xi_{\lambda},
\dot w^{(l)} \cdot \xi_{\lambda} \right> \\
& = & \left< x_{i_l}(m_l)\dot s_{i_l}^{-1}\dot w^{(l+1)}\cdot \xi_{\lambda},
\dot w^{(l)} \cdot \xi_{\lambda} \right>
\left< g^{(l+1)}\cdot \xi_{\lambda},\dot w^{(l+1)}\cdot \xi_{\lambda}
\right> \\
& = & \left< \dot s_{i_l}^{-1}\dot w^{(l+1)}\cdot \xi_{\lambda},
\dot w^{(l)} \cdot \xi_{\lambda} \right>
\left< g^{(l+1)}\cdot \xi_{\lambda},\dot w^{(l+1)}\cdot \xi_{\lambda}
\right> \\
& = &
(-1)^{\left< s_{i_{l+1}}\cdots s_{i_n}\lambda,\alpha^{\vee}_{i_l}\right>}
\Delta_{\lambda}^{w^{(l+1)}\lambda}(g_{(l+1)}).
\end{eqnarray*}
The last equality follows from the fact that
$\dot s_{i_l}^{-1}=\alpha_{i_l}^{\vee}(-1)\dot s_{i_l}$.

The result for $l$ now follows (in each case) from the inductive
hypothesis and we are done.

(2) We have:
\begin{eqnarray*}
\Delta_{\omega_{i_k}}^{v_{(k-1)}\omega_{i_k}}(g_{(k)}) & = &
\left< g_{(k)} \cdot \xi_{\omega_{i_k}},\dot v_{(k-1)}\cdot \xi_{\omega_{i_k}}
\right> \\
& = & \left< g_{(k-1)}x_{i_k}(m_k)\dot s_{i_k}^{-1}
\cdot \xi_{\omega_{i_k}},\dot v_{(k-1)}\cdot \xi_{\omega_{i_k}}\right> \\
& = & -\left< g_{(k-1)}\dot s_{i_k}\cdot \xi_{\omega_{i_k}},
\dot v_{(k-1)}\cdot \xi_{\omega_{i_k}}\right> -
m_k\left< g_{(k-1)}\cdot \xi_{\omega_{i_k}},\dot v_{(k-1)}
\cdot \xi_{\omega_{i_k}}\right>, \\
& = & -\Delta^{v_{(k-1)}\omega_{i_k}}_{s_{i_k}\omega_{i_k}}(g_{(k-1)})-m_k,
\end{eqnarray*}
noting that, since $g_{(k-1)}\in U^-\dot v_{(k-1)}$, we have that
$\langle g_{(k-1)}\cdot \xi_{\omega_{i_k}}, \dot v_{(k-1)}
\cdot \xi_{\omega_{i_k}}\rangle=1$. The result follows.
\end{proof}

\begin{rems} \label{r:formulas}
(1) Let $t_k=-1$ for $k\in J^-_{\mathbf{v}}$ and let $t_k=1$ for
$k\in J^+_{\mathbf{v}}$ (so that now
$t_k$ is defined for $k=1,2,\ldots ,n$). Then the formula in
Lemma~\ref{l:gminor}(1) can be rewritten as:

$$\Delta_{\lambda}^{w_{(k)}\lambda}(g_{(k)})=
\prod_{l=1}^k
t_l^{\left< s_{i_{l+1}}s_{i_{l+2}}\cdots s_{i_k}\lambda,\alpha^{\vee}_{i_l}
\right>}.$$

(2) The following lemma gives an expression for $m_k$ which is simpler
than the Chamber Ansatz version, Theorem~\ref{t:ansatz}(2). However,
the Chamber Ansatz formula for the $m_k$ will be more useful in
Section~\ref{s:coordinates}.
\end{rems}

\begin{lem}
For $k\in J^-_{\mathbf{v}}$, we have:
$$m_k=-\frac{
\Delta_{w_{(k)}\omega_{i_k}}^{v_{(k-1)}\omega_{i_k}}(z)}
{\Delta_{w_{(k)}\omega_{i_k}}^{v_{(k)}\omega_{i_k}}(z)}
-\Delta^{v_{(k-1)}\omega_{i_k}}_{s_{i_k}\omega_{i_k}}(g_{(k-1)}).$$
\end{lem}
\begin{proof}
This follows immediately from Lemma~\ref{l:crossminor} and
Lemma~\ref{l:gminor}(2).
\end{proof}
\begin{proof}[Proof of Theorem~\ref{t:ansatz}]
We can now prove Theorem~\ref{t:ansatz} by using Lemmas~\ref{l:crossminor}
and~\ref{l:gminor} and Remark~\ref{r:formulas}(1) to substitute for the minors
appearing in the expressions on the right hand sides of
Theorem~\ref{t:ansatz}(1) and (2).
We first claim that, for $k=1,2,\ldots ,n$, we have:
\begin{equation} \label{e:generalansatz}
t_k=
\frac{
\prod_{j\not=i_k}\Delta_{w_{(k)}\omega_j}^{v_{(k)}\omega_j}(z)^{-a_{j,i_k}}}
{\Delta_{w_{(k)}\omega_{i_k}}^{v_{(k)}\omega_{i_k}}(z)\Delta_{w_{(k-1)}\omega_{i_k}}^{v_{(k-1)}\omega_{i_k}}(z)}.
\end{equation}

(See~\cite[4.3]{BeZel:TotPos} for a similar proof of this statement in the
special case where $J^+_{\mathbf{v}}=\emptyset$).
We have, using Lemma~\ref{l:crossminor}(1) and Remark~\ref{r:formulas}(1):
\begin{equation*}
\frac{
\prod_{j\not=i_k}\Delta_{w_{(k)}\omega_j}^{v_{(k)}\omega_j}(z)^{-a_{j,i_k}}}
{\Delta_{w_{(k)}\omega_{i_k}}^{v_{(k)}\omega_{i_k}}(z)\Delta_{w_{(k-1)}
\omega_{i_k}}^{v_{(k-1)}\omega_{i_k}}(z)}
=
\frac{\prod_{j\not=i_k} \prod_{l=1}^k t_l^{a_{j,i_k}\left< s_{i_{l+1}}\cdots
s_{i_k}\omega_j,\alpha_{i_l}^{\vee} \right> }}
{\prod_{l=1}^k t_l^{-\left< s_{i_{l+1}}\cdots s_{i_k}\omega_{i_k},
\alpha_{i_l}^{\vee}\right> }
\prod_{l=1}^{k-1}t_l^{-\left< s_{i_{l+1}}\cdots s_{i_{k-1}}\omega_{i_k},
\alpha_{i_l}^{\vee}\right> }}.
\end{equation*}
The exponent of $t_k$ is given by
\begin{eqnarray*}
\sum_{j\not=i_k}a_{j,i_k}\left< \omega_j,\alpha_{i_k}^{\vee} \right>
+\left< \omega_{i_k},\alpha_{i_k}^{\vee}\right>
& = &
\left<\sum_{j\not=i_k}a_{j,i_k}\omega_j+\omega_{i_k},\alpha_{i_k}^{\vee}\right>
\\
& = & \left< \alpha_{i_k}-\omega_{i_k},\alpha_{i_k}^{\vee}\right>=1.
\end{eqnarray*}
If $k<k'$, the exponent of $t_{k'}$ is clearly zero. If $k'<k$, then
the exponent of $t_{k'}$ is given by
\begin{equation*}
\left( \sum_{j\not=i_k}a_{j,i_k}\left< s_{i_{k'+1}}\cdots s_{i_k}\omega_j,
\alpha_{i_{k'}}^{\vee} \right>\right)
+\left< s_{i_{k'+1}}\cdots s_{i_k}\omega_{i_k},\alpha_{i_{k'}}^{\vee}\right>
+\left< s_{i_{k'+1}}\cdots
s_{i_{k-1}}\omega_{i_k},\alpha_{i_{k'}}^{\vee}\right>
\end{equation*}
\begin{eqnarray*}
&=&\left< s_{i_{k'+1}}\cdots s_{i_k}\left(
\left(\sum_{j\not=i_k}a_{j,i_k}\omega_j\right)+2\omega_{i_k}-\alpha_{i_k}\right),
\alpha_{i_{k'}}^{\vee} \right>
\\
&=&\left< s_{i_{k'+1}}\cdots s_{i_k}(
\alpha_{i_k}-\alpha_{i_k}),\alpha_{i_{k'}}^{\vee} \right>=0,
\end{eqnarray*}
and the claim~(\ref{e:generalansatz}) is proved; Theorem~\ref{t:ansatz}(1)
is a special case.

We now prove Theorem~\ref{t:ansatz}(2). Suppose that $k\in J^-_{\mathbf{v}}$.
Using Lemma~\ref{l:crossminor} and~(\ref{e:generalansatz})
(noting that $t_k=-1$), we see that
\begin{eqnarray*}
\frac{
\Delta_{w_{(k)}\omega_{i_k}}^{v_{(k-1)}\omega_{i_k}}(z)
\Delta_{w_{(k-1)}\omega_{i_k}}^{v_{(k-1)}\omega_{i_k}}(z)}
{\prod_{j\not=i_k}\Delta_{w_{(k)}\omega_j}^{v_{(k)}\omega_j}(z)^{-a_{j,i_k}}}
& = &
\frac{\Delta_{\omega_{i_k}}^{v_{(k-1)}\omega_{i_k}}(g_{(k)})
\Delta_{w_{(k-1)}\omega_{i_k}}^{v_{(k-1)}\omega_{i_k}}(z)}
{\Delta_{\omega_{i_k}}^{w_{(k)}\omega_{i_k}}(g_{(k)})
\prod_{j\not=i_k}\Delta_{w_{(k)}\omega_j}^{v_{(k)}\omega_j}(z)^{-a_{j,i_k}}}
\\
& = &
\frac{\Delta_{\omega_{i_k}}^{v_{(k-1)}\omega_{i_k}}(g_{(k)})
\Delta_{w_{(k)}\omega_{i_k}}^{v_{(k)}\omega_{i_k}}(z)
\Delta_{w_{(k-1)}\omega_{i_k}}^{v_{(k-1)}\omega_{i_k}}(z)}
{\prod_{j\not=i_k}\Delta_{w_{(k)}\omega_j}^{v_{(k)}\omega_j}(z)^{-a_{j,i_k}}}
\\
& = &
-\Delta_{\omega_{i_k}}^{v_{(k-1)}\omega_{i_k}}(g_{(k)}).
\end{eqnarray*}
Theorem~\ref{t:ansatz}(2) now follows from Lemma~\ref{l:gminor}(2),
and the proof of Theorem~\ref{t:ansatz} is complete.
\end{proof}
\section{A change of coordinates}\label{s:coordinates}
We can gain some more insight into the structure of the formulas
from Theorem~\ref{t:ansatz} if we consider the standard and
special chamber minors as providing an alternative system
of coordinates on $\mathcal R_{\v,\w}$.
\begin{prop}\label{p:ansatz} Let $\v\prec \w$. With notation as
above the map
 \begin{align*}
 \mathcal R_{\v,\w}&\To (\k^*)^{J^\circ_\v}\x (\k)^{J^-_\v} \\
 z\dot w\cdot B^+ & \mapsto
\left( (\Delta^{v_{(j)}\omega_{i_j}}_{w_{(j)}\omega_{i_j}}(z)
)_{j\in J^\circ_\v },
(\Delta^{v_{(j-1)}\omega_{i_j}}_{w_{(j)}\omega_{i_j}}(z) )_{j\in
J^-_\v} \right)
 \end{align*}
is an isomorphism.
\end{prop}
For the special case $\v=(1,\dotsc,1)$ see also
Theorem~4.3 and Corollary~4.4 in \cite{BeZel:TotPos}.
\begin{proof}
Let $\w_{(k)}=(w_{(0)},\dotsc, w_{(k)})$ be the reduced expression
for $w_{(k)}$ obtained from $\w$ by
truncation, and $\v_{(k)}$ the corresponding truncation of $\v$.
The proof of the proposition is by induction on $k$ and using
Theorem~\ref{t:ansatz}. The start of the induction is trivial, so
let us assume the proposition is true for
$\mathcal R_{\v_{(k-1)},\w_{(k-1)}}$. We have three cases for $k$.
\begin{enumerate}
\item
If $k\in J^\circ_\v$, then as in the proof of
Proposition~\ref{p:parameterization} we have
  \begin{equation*}\label{e:iso1}
  \begin{array}{lcl}
 \mathcal R_{\v_{(k-1)},\w_{(k-1)}}\x \k^*  & \overset \sim \To &
 \mathcal R_{\v_{(k)},\w_{(k)}},
 \\
 (g_{(k-1)}\cdot B^+, t_k) & \mapsto & g_{(k-1)}y_{i_k}(t_k)
 \cdot B^+, \qquad \quad g_{(k-1)}\in G_{\v_{(k-1)},\w_{(k-1)}} .
 \end{array}
  \end{equation*}
Compose this map with
 \begin{equation}\label{e:map1}
  \begin{array}{lcl}
 \mathcal R_{\v_{(k)},\w_{(k)}}  &\To &
 \mathcal R_{\v_{(k-1)},\w_{(k-1)}} \x \k^*,    \\
 z\dot w_{(k)}\cdot B^+ &\mapsto &
 \left(z\dot w_{(k-1)}\cdot B^+,
 \Delta^{v_{(k)}\omega_{i_k}}_{w_{(k)}\omega_{i_k}}(z)\right )
 \end{array}
  \end{equation}
to get a map
$\psi_k:\mathcal R_{\v_{(k-1)},\w_{(k-1)}}\x \k^*\to
\mathcal R_{\v_{(k-1)},\w_{(k-1)}}\x \k^*$.
The Chamber Ansatz says that $t_k$ can be recovered
from $z\dot w_{(k)}\cdot B^+$ by
 \begin{equation}\label{e:tk2}
 t_k=
 a_k\left (z\dot w_{(k-1)}\cdot B^+\right )\
 \Delta^{v_{(k)}\omega_{i_k}}_{w_{(k)}\omega_{i_k}}(z)\inv,
 \end{equation}
where
 \begin{equation*}\label{e:ak1}
a_k(z\dot w_{(k-1)}\cdot B^+):= \frac{
\prod_{j\not=i_k}\Delta_{w_{(k)}\omega_j}^{v_{(k)}\omega_j}(z)^{-a_{j,i_k}}
}{\Delta_{w_{(k-1)}\omega_{i_k}}^{v_{(k-1)}\omega_{i_k}}(z)},
 \qquad z\dot w_{(k-1)}\cdot B^+\in\mathcal R_{\v_{(k-1)},\w_{(k-1)}}.
 \end{equation*}
Note that this gives a well-defined map
$a_k:\mathcal R_{\v_{(k-1)},\w_{(k-1)}}\to \k^*$, since $a_k$
is made up of standard chamber minors for $(\v_{(k-1)},\w_{(k-1)})$.
Now the formula \eqref{e:tk2} gives rise to an inverse to $\psi_k$ .
Hence also \eqref{e:map1} is
an isomorphism, and the proposition holds for $\mathcal R_{\v_{(k)},\w_{(k)}}$
by the induction hypothesis.

\item
Suppose $k\in J^-_\v$. Then we have an isomorphism
  \begin{equation*}\label{e:iso2}
  \begin{array}{lcl}
 \mathcal R_{\v_{(k-1)},\w_{(k-1)}}\x \k  & \overset \sim \To &
 \mathcal R_{\v_{(k)},\w_{(k)}},
 \\
 (g_{(k-1)}\cdot B^+, m_k) & \mapsto & g_{(k-1)}x_{i_k}(m_k)\dot s_{i_k}\inv
 \cdot B^+, \qquad \quad g_{(k-1)}\in G_{\v_{(k-1)},\w_{(k-1)}},
 \end{array}
  \end{equation*}
which we can compose with
 \begin{equation}\label{e:map2}
  \begin{array}{lcl}
 \mathcal R_{\v_{(k)},\w_{(k)}}  &\To &
 \mathcal R_{\v_{(k-1)},\w_{(k-1)}} \x \k,    \\
 z\dot w_{(k)}\cdot B^+ &\mapsto &
 \left(z\dot w_{(k-1)}\cdot B^+,
 \Delta^{v_{(k-1)}\omega_{i_k}}_{w_{(k)}\omega_{i_k}}(z)\right )
 \end{array}
  \end{equation}
to get a map  $\psi_k:\mathcal R_{\v_{(k-1)},\w_{(k-1)}}\x \k\to
\mathcal R_{\v_{(k-1)},\w_{(k-1)}}\x \k$.
By Theorem~\ref{t:ansatz} one can recover the $m_k$ coordinate from
$z\dot w_{(k)}\cdot B^+$ by
 \begin{equation}\label{e:mkequation}
m_k=a_k\left(z\dot w_{(k-1)}\cdot B^+\right)
\Delta^{v_{(k-1)}\omega_{i_k}}_{w_{(k)}\omega_{i_k}}(z) -
 b_k\left(z\dot w_{(k-1)}\cdot B^+\right),
 \end{equation}
where $a_k:\mathcal R_{\v_{(k-1)},\w_{(k-1)}}\to\k^*$
and $b_k:\mathcal R_{\v_{(k-1)},\w_{(k-1)}}\to\k$ are given by
\begin{equation*}
\begin{array}{lcll}
 a_k(z\dot w_{(k-1)}\cdot B^+)&=& \frac{
 \Delta_{w_{(k-1)}\omega_{i_k}}^{v_{(k-1)}\omega_{i_k}}(z)}
{\prod_{j\not=i_k}\Delta_{w_{(k)}\omega_j}^{v_{(k)}\omega_j}(z)^{-a_{j,i_k}}},
 &\qquad\qquad z\dot w_{(k-1)}\cdot B^+\in\mathcal R_{\v_{(k-1)},\w_{(k-1)}},
\\
 b_k(g_{(k-1)}\cdot B^+)&=&
\Delta_{s_{i_k}\omega_{i_k}}^{v_{(k-1)}\omega_{i_k}}(g_{(k-1)}),
 &\qquad\qquad g_{(k-1)}\in G_{\v_{(k-1)},\w_{(k-1)}}.
\end{array}
\end{equation*}
Now the identity \eqref{e:mkequation} gives an inverse to the map
$\psi_k$. So \eqref{e:map2} is an isomorphism and the
proposition holds for $\mathcal R_{\v_{(k)},\w_{(k)}}$.
\item
If $k\in J^+_\v$ then $\mathcal R_{\v_{(k)},\w_{(k)}}\cong
\mathcal R_{\v_{(k-1)},\w_{(k-1)}}$ and we are done.
\end{enumerate}
\end{proof}
\begin{rem}\label{r:ansatz}
Assuming Proposition~\ref{p:ansatz},
much of the structure of Theorem~\ref{t:ansatz}
is already determined. The proposition says that we may take the
$(\Delta^{v_{(j)}\omega_{i_j}}_{w_{(j)}\omega_{i_j}}(z) )_{j\in
J^\circ_\v }$ in $\k^*$, and the
$(\Delta^{v_{(j-1)}\omega_{i_j}}_{w_{(j)}\omega_{i_j}}(z) )_{j\in
J^-_\v}$ in $\k$ together as coordinates for $\mathcal R_{\v,\w}$.
And then Theorem~\ref{t:ansatz} can be interpreted roughly as the
transition from these coordinates to the coordinates $(t_j)_{j\in
J^\circ_\v }$ and $(m_j)_{j\in J^-_\v }$ from the factorization.
>From the outset these two sets of coordinates have to be quite
closely related, since they
are both compatible with reduction. In either setting
the map
$\pi^{w}_{w_{(k)}}:\mathcal R_{\v,\w}\to\mathcal R_{\v_{(k)},\w_{(k)}}$
corresponds to the projection onto the first $k$ coordinates.

Let us consider explicitly an element $g\cdot B^+=z\dot w\cdot  B^+$
with fixed reduction
$g_{(k)}\cdot B^+=z\dot w_{(k)}\cdot B^+$. Then this is
equivalent to fixing coordinates
$t_j$ and $m_j$, or to fixing the minors
$\Delta^{v_{(j)}\omega_{i_j}}_{w_{(j)}\omega_{i_j}}(z)$ and
$\Delta^{v_{(j-1)}\omega_{i_j}}_{w_{(j)}\omega_{i_j}}(z)$, where
$j\le k$ in $J^\circ_\v$ or $J^-_\v$.
If $k+1\in J^\circ_\v$ then the change of coordinate from
$\Delta^{v_{(k+1)}\omega_{i_{k+1}}}_{w_{(k+1)}\omega_{i_{k+1}}}(z)$
to $t_{k+1}$ amounts to an invertible map $\k^*\to \k^*$, which depends
only on the earlier coordinates. So it has to be of the form
$z\mapsto a z^{\pm 1}$ for some $a\in \k^*$. The Chamber Ansatz simply
says the map is of the form $z\mapsto a z\inv$ and describes the
coefficient $a$ explicitly in terms of earlier chamber minors of $z$.

If $k\in J^-_\v$ then the change of the coordinate
$\Delta^{v_{(k)}\omega_{i_{k+1}}}_{w_{(k+1)}\omega_{i_{k+1}}}(z)$
to $m_{k+1}$ amounts to an invertible map $\k\to \k$, which depends
only on the earlier coordinates. Therefore this map must be of the form
$z\mapsto a z +b$ for $a\in \k^*$ and $b\in \k$.
Here again $a$ is computed by the Chamber Ansatz,
and $b$ is the correction term in Theorem~\ref{t:ansatz}.
\end{rem}
\section{The generalized Chamber Ansatz for $SL_d$} \label{s:ansatztypea}
Suppose we are given $B=z\dot w\cdot B^+$, with
$z\in U^+$, with fixed reduced expression $\mathbf{w}$ for $w$.
We can use Proposition~\ref{p:algorithm} to determine
which Deodhar component ${\mathcal R}_{\mathbf{v},\mathbf{w}}$ contains $B$,
where
$\mathbf{v}$ is a distinguished subexpression for $v$ in $\mathbf{w}$.
We shall give an explicit example in Section~\ref{s:componentexample} of how
to do this.
Then $B=g\cdot B^+$, where $g=g_1g_2\cdots g_n$, and
$$g_k=\left\{\begin{array}{ll}
y_{i_k}(t_k) & k\in J^{\circ}_{\mathbf{v}}, \\
\dot s_{i_k} & k\in J^+_{\mathbf{v}}, \\
x_{i_k}(m_k) \dot s_{i_k}^{-1} & k\in J^-_{\mathbf{v}}.
\end{array}\right.
$$
The generalized Chamber Ansatz (Theorem~\ref{t:ansatz}) gives formulas for
the $t_k$ and $m_k$ in terms of minors of $z$ (and the minor
$\Delta^{v_{(k-1)}\omega_{i_k}}_{s_{i_k}\omega_{i_k}}(g_{(k-1)})$ of
$g_{(k-1)}$). We write Theorem~\ref{t:ansatz}(2) in the form
$$m_k=r_k-\Delta^{v_{(k-1)}\omega_{i_k}}_{s_{i_k}\omega_{i_k}}(g_{(k-1)}),$$
where
$$
r_k=\frac{
\Delta_{w_{(k)}\omega_{i_k}}^{v_{(k-1)}\omega_{i_k}}(z)\Delta_{w_{(k-1)}\omega_{i_k}}^{v_{(k-1)}\omega_{i_k}}(z)}
{\prod_{j\not=i_k}\Delta_{w_{(k)}\omega_j}^{v_{(k)}\omega_j}(z)^{-a_{j,i_k}}}.
$$
This is the term $a_k\left(z\dot w_{(k-1)}\cdot B^+\right)
\Delta^{v_{(k-1)}\omega_{i_k}}_{w_{(k)}\omega_{i_k}}(z)$
in equation~(\ref{e:mkequation}).
In this section, we give a graphical algorithm (generalizing that
of~\cite{BeFoZe:TotPos}) for determining the coefficients $t_k$ and $r_k$,
in the case where $G=SL_d$. The coefficients $m_k$ can then be computed
by computing the minors
$\Delta^{v_{(k-1)}\omega_{i_k}}_{s_{i_k}\omega_{i_k}}(g_{(k-1)})$
inductively, noting that $g_{(k-1)}$ depends only on the coefficients
$t_j$ and $m_j$ for $j\leq k-1$ (an example of this will be given in
section~\ref{s:example}).

We employ a generalised version of the pseudoline arrangements used
in~\cite{BeFoZe:TotPos}, in which two pseudolines can
either intersect, as in~\cite{BeFoZe:TotPos}, or pass over or under each other
(see below for examples). These can also be regarded as diagrams of singular
braids~\cite{Baez:Link,Birman:NewView}.

\begin{figure}[htbp]
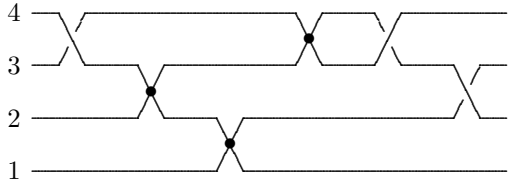

\beginpicture

\setcoordinatesystem units <0.35cm,0.35cm>             
\setplotarea x from -14 to 18, y from 0 to 8       
\linethickness=0.5pt           

\setlinear \plot 0 0  7 0  / %
\setlinear \plot 7 0  8 2  / %
\setlinear \plot 8 2  16 2  / %
\setlinear \plot 16 2 16.35 2.7  / %
\setlinear \plot 16.65 3.3 17 4  / %
\setlinear \plot 17 4  18 4  / %

\setlinear \plot 0 2  4 2  / %
\setlinear \plot 4 2  5 4  / %
\setlinear \plot 5 4  10 4 / %
\setlinear \plot 10 4 11 6  / %
\setlinear \plot 11 6 13 6  / %
\setlinear \plot 13 6 13.35 5.3  / %
\setlinear \plot 13.65 4.7  14 4  / %
\setlinear \plot 14 4  16 4  / %
\setlinear \plot 16 4  17 2  / %
\setlinear \plot 17 2  18 2  / %

\setlinear \plot 0 4  1 4  / %
\setlinear \plot 1 4  1.35 4.7  / %
\setlinear \plot 1.65 5.3  2 6  / %
\setlinear \plot 2 6  10 6  / %
\setlinear \plot 10 6 11 4  / %
\setlinear \plot 11 4 13 4  / %
\setlinear \plot 13 4 14 6  / %
\setlinear \plot 14 6 18 6  / %

\setlinear \plot 0 6  1 6  / %
\setlinear \plot 1 6  2 4  / %
\setlinear \plot 2 4  4 4  / %
\setlinear \plot 4 4  5 2  / %
\setlinear \plot 5 2  7 2  / %
\setlinear \plot 7 2  8 0  / %
\setlinear \plot 8 0  18 0  / %

\put{$\bullet$}[c] at 4.5 3
\put{$\bullet$}[c] at 7.5 1
\put{$\bullet$}[c] at 10.5 5

\put{$1$}[c] at -0.7 0
\put{$2$}[c] at -0.7 2
\put{$3$}[c] at -0.7 4
\put{$4$}[c] at -0.7 6

\endpicture
\caption{Ansatz arrangement (unlabeled) for
$\underline{s_3} s_2 s_1 \underline{s_3 s_2}$.
Note that $g=\dot s_3 y_2(t_2) y_1(t_3) x_3(m_4) \dot s_3^{-1} \dot s_2$.}
\label{f:emptyansatzarrangement}
\end{figure}

The main idea is to associate such an arrangement (which we call the
{\em ansatz arrangement}) to the pair
$\mathbf{v},\mathbf{w}$. For example, if
$\mathbf{w}=(1,s_3,s_3s_2,s_3s_2s_1,s_3s_2s_1s_3,s_3s_2s_1s_3s_2)$, and
$\mathbf{v}$ is the distinguished subexpression $(1,s_3,s_3,s_3,1,s_2)$
for $s_2$ in $\mathbf{w}$, then the arrangement is as in
Figure~\ref{f:emptyansatzarrangement}.

The pair $\mathbf{v}$, $\mathbf{w}$ determines the factors $g_k$ of
$g$, which in turn determine the ansatz arrangement in the
following way.
It consists of $d$ pseudolines, numbered $1,2,\ldots ,d$, from
bottom to top on the left hand side of the arrangement. Each factor
$x_{i_k}(m_k)$, $y_{i_k}(t_k)$, $\dot s_{i_k}$ or $\dot s_{i_k}^{-1}$ of
$g$ gives rise to a constituent of the arrangement, in which pseudolines
at level $i_k$ from the bottom of the arrangement are either braided
or cross at singular point.
Note that $x_{i_k}(m_k)$ and $\dot s_{i_k}^{-1}$ are treated as
separate factors. The rules for how this is done are given in
Figure~\ref{f:ansatzrules}.

\begin{figure}[htbp]
\beginpicture

\setcoordinatesystem units <0.35cm,0.35cm>             
\setplotarea x from -8.5 to 33, y from 6 to 14       
\linethickness=0.2pt           


\put{Factor of $g$}[c] at 1 13
\put{Constituent of}[c] at 1 10
\put{ansatz}[c] at 1 9
\put{arrangement}[c] at 1 8

\put{$x_{i_k}(m_k)$}[c] at 8.5 13
\put{$y_{i_k}(t_k)$}[c] at 15.5 13
\put{$\dot s_{i_k}$}[c] at 22.5 13
\put{$\dot s_{i_k}^{-1}$}[c] at 29.5 13

\setlinear \plot -2.5 12  33 12 / %
\setlinear \plot 5 6   5 14 / %
\setlinear \plot 12 6  12 14 / %
\setlinear \plot 19 6  19 14 / %
\setlinear \plot 26 6  26 14 / %

\linethickness=0.7pt           

\setlinear \plot 7 8  8 8 / %
\setlinear \plot 8 8  9 10 / %
\setlinear \plot 9 10  10 10 / %
\setlinear \plot 7 10  8 10 / %
\setlinear \plot 8 10  9 8 / %
\setlinear \plot 9 8  10 8 / %
\put{$\bullet$}[c] at 8.5 9

\setlinear \plot 14 8  15 8 / %
\setlinear \plot 15 8  16 10 / %
\setlinear \plot 16 10  17 10 / %
\setlinear \plot 14 10  15 10 / %
\setlinear \plot 15 10  16 8 / %
\setlinear \plot 16 8  17 8 / %
\put{$\bullet$}[c] at 15.5 9

\setlinear \plot 21 8  22 8 / %
\setlinear \plot 22 8  22.35 8.7 / %
\setlinear \plot 22.65 9.3  23 10 / %
\setlinear \plot 23 10  24 10 / %
\setlinear \plot 21 10  22 10 / %
\setlinear \plot 22 10  23 8 / %
\setlinear \plot 23 8  24 8 / %

\setlinear \plot 28 8  29 8 / %
\setlinear \plot 29 8  30 10 / %
\setlinear \plot 30 10  31 10 / %
\setlinear \plot 28 10  29 10 / %
\setlinear \plot 29 10  29.35 9.3 / %
\setlinear \plot 29.65 8.7  30 8 / %
\setlinear \plot 30 8  31 8 / %

\endpicture
\caption{Constituents of the ansatz arrangement.}
\label{f:ansatzrules}
\end{figure}

\begin{figure}[htbp]
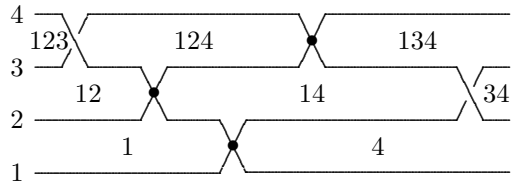

\beginpicture

\setcoordinatesystem units <0.35cm,0.35cm>             
\setplotarea x from -14 to 18, y from 0 to 8       
\linethickness=0.5pt           

\setlinear \plot 0 0  7 0  / %
\setlinear \plot 7 0  8 2  / %
\setlinear \plot 8 2  16 2  / %
\setlinear \plot 16 2 16.35 2.7  / %
\setlinear \plot 16.65 3.3 17 4  / %
\setlinear \plot 17 4  18 4  / %

\setlinear \plot 0 2  4 2  / %
\setlinear \plot 4 2  5 4  / %
\setlinear \plot 5 4  10 4 / %
\setlinear \plot 10 4 11 6  / %
\setlinear \plot 11 6 18 6  / %

\setlinear \plot 0 4  1 4  / %
\setlinear \plot 1 4  1.35 4.7  / %
\setlinear \plot 1.65 5.3  2 6  / %
\setlinear \plot 2 6  10 6  / %
\setlinear \plot 10 6 11 4  / %
\setlinear \plot 11 4 16 4  / %
\setlinear \plot 16 4 17 2  / %
\setlinear \plot 17 2 18 2  / %

\setlinear \plot 0 6  1 6  / %
\setlinear \plot 1 6  2 4  / %
\setlinear \plot 2 4  4 4  / %
\setlinear \plot 4 4  5 2  / %
\setlinear \plot 5 2  7 2  / %
\setlinear \plot 7 2  8 0  / %
\setlinear \plot 8 0  18 0  / %

\put{$\bullet$}[c] at 4.5 3
\put{$\bullet$}[c] at 7.5 1
\put{$\bullet$}[c] at 10.5 5

\put{$123$}[c] at 0.5 5
\put{$124$}[c] at 6 5
\put{$134$}[c] at 14.5 5
\put{$12$}[c] at 2 3
\put{$14$}[c] at 10.5 3
\put{$34$}[c] at 17.5 3
\put{$1$}[c] at 3.5 1
\put{$4$}[c] at 13 1

\put{$1$}[c] at -0.7 0
\put{$2$}[c] at -0.7 2
\put{$3$}[c] at -0.7 4
\put{$4$}[c] at -0.7 6

\endpicture
\caption{Lower arrangement for
$\underline{s_3} s_2 s_1 \underline{s_3 s_2}$.
Note that $g=\dot s_3 y_2(t_2) y_1(t_3) x_3(m_4) \dot s_3^{-1} \dot s_2$.}
\label{f:examplelowerarrangement}
\end{figure}

As usual, we define a {\em chamber} of a generalised pseudoline arrangement
to be a component of the complement of the union of the pseudolines in the
arrangement (for this definition we interpret the under- and over-crossings
as singular points). In order to label the chambers, we need two
auxilliary pseudoline arrangements associated to the pair
$\mathbf{v},\mathbf{w}$, which we call the {\em upper} and
{\em lower arrangements} (since, as will be seen, they will determine
upper and lower subscripts of chamber minors).
These arrangements are defined in the same way as the ansatz arrangement,
except that different rules for the factors of $g$ are employed. These
are described in Figure~\ref{f:upperlowerrules}, and the upper and
lower arrangements for the example above are given in
Figures~\ref{f:exampleupperarrangement} and~\ref{f:examplelowerarrangement}.
The chambers for these arrangements are labeled with the labels of the
strings passing below them.

\begin{figure}[htbp]
\beginpicture

\setcoordinatesystem units <0.35cm,0.35cm>             
\setplotarea x from -8 to 33, y from 0 to 14       
\linethickness=0.2pt           


\put{Factor of $g$}[c] at 1 13
\put{Constituent of}[c] at 1 10
\put{upper}[c] at 1 9
\put{arrangement}[c] at 1 8
\put{Constituent}[c] at 1 4
\put{lower}[c] at 1 3
\put{arrangement}[c] at 1 2

\put{$x_{i_k}(m_k)$}[c] at 8.5 13
\put{$y_{i_k}(t_k)$}[c] at 15.5 13
\put{$\dot s_{i_k}$}[c] at 22.5 13
\put{$\dot s_{i_k}^{-1}$}[c] at 29.5 13

\setlinear \plot -2.5 12  33 12 / %
\setlinear \plot -2.5 6   33 6  / %
\setlinear \plot 5 0   5 14 / %
\setlinear \plot 12 0  12 14 / %
\setlinear \plot 19 0  19 14 / %
\setlinear \plot 26 0  26 14 / %

\linethickness=0.7pt           

\setlinear \plot 7 2  8 2 / %
\setlinear \plot 8 2  9 4 / %
\setlinear \plot 9 4  10 4 / %
\setlinear \plot 7 4  8 4 / %
\setlinear \plot 8 4  9 2 / %
\setlinear \plot 9 2  10 2 / %
\put{$\bullet$}[c] at 8.5 3

\setlinear \plot 7 8  10 8 / %
\setlinear \plot 7 10  10 10 / %

\setlinear \plot 14 2  15 2 / %
\setlinear \plot 15 2  16 4 / %
\setlinear \plot 16 4  17 4 / %
\setlinear \plot 14 4  15 4 / %
\setlinear \plot 15 4  16 2 / %
\setlinear \plot 16 2  17 2 / %
\put{$\bullet$}[c] at 15.5 3

\setlinear \plot 14 8  17 8 / %
\setlinear \plot 14 10  17 10 / %

\setlinear \plot 21 2  22 2 / %
\setlinear \plot 22 2  22.35 2.7 / %
\setlinear \plot 22.65 3.3  23 4 / %
\setlinear \plot 23 4  24 4 / %
\setlinear \plot 21 4  22 4 / %
\setlinear \plot 22 4  23 2 / %
\setlinear \plot 23 2  24 2 / %

\setlinear \plot 21 8  22 8 / %
\setlinear \plot 22 8  22.35 8.7 / %
\setlinear \plot 22.65 9.3  23 10 / %
\setlinear \plot 23 10  24 10 / %
\setlinear \plot 21 10  22 10 / %
\setlinear \plot 22 10  23 8 / %
\setlinear \plot 23 8  24 8 / %

\setlinear \plot 28 2  31 2 / %
\setlinear \plot 28 4  31 4 / %

\setlinear \plot 28 8  29 8 / %
\setlinear \plot 29 8  30 10 / %
\setlinear \plot 30 10  31 10 / %
\setlinear \plot 28 10  29 10 / %
\setlinear \plot 29 10  29.35 9.3 / %
\setlinear \plot 29.65 8.7  30 8 / %
\setlinear \plot 30 8  31 8 / %

\endpicture
\caption{Constituents of the auxilliary arrangements.}
\label{f:upperlowerrules}
\end{figure}

\begin{figure}[htbp]
\beginpicture

\setcoordinatesystem units <0.35cm,0.35cm>             
\setplotarea x from -14 to 18, y from 0 to 8       
\linethickness=0.5pt           

\setlinear \plot  0 0  18 0 / %

\setlinear \plot  0 2  16 2 / %
\setlinear \plot  16 2  16.35 2.7 / %
\setlinear \plot  16.65 3.3  17 4 / %
\setlinear \plot  17 4  18 4 / %

\setlinear \plot 0 4  1 4  / %
\setlinear \plot 1 4  1.35 4.7  / %
\setlinear \plot 1.65 5.3  2 6 / %
\setlinear \plot 2 6  13 6  / %
\setlinear \plot 13 6 13.35 5.3  / %
\setlinear \plot 13.65 4.7  14 4  / %
\setlinear \plot 14 4  16 4  / %
\setlinear \plot 16 4  17 2  / %
\setlinear \plot 17 2  18 2  / %

\setlinear \plot 0 6  1 6 / %
\setlinear \plot 1 6  2 4 / %
\setlinear \plot 2 4  13 4 / %
\setlinear \plot 13 4 14 6 / %
\setlinear \plot 14 6 18 6 / %

\put{$123$}[c] at 0.5 5
\put{$124$}[c] at 7.5 5
\put{$123$}[c] at 16  5
\put{$12$}[c] at 8 3
\put{$13$}[c] at 17.5 3
\put{$1$}[c] at 9 1

\put{$1$}[c] at -0.7 0
\put{$2$}[c] at -0.7 2
\put{$3$}[c] at -0.7 4
\put{$4$}[c] at -0.7 6

\endpicture
\caption{Upper arrangement for
$\underline{s_3} s_2 s_1 \underline{s_3 s_2}$.
Note that $g=\dot s_3 y_2(t_2) y_1(t_3) x_3(m_4) \dot s_3^{-1} \dot s_2$.}
\label{f:exampleupperarrangement}
\end{figure}

We note that, since $G=SL_d$, the generalized minors of
Definition~\ref{d:minors} coincide with the usual minors of matrices.
We denote by $\Delta^R_S$ the minor with row set $R$ and column set $S$
(interpreted as $1$ if $R=S=\emptyset$).
Suppose that $X$ is a chamber of the ansatz arrangement. Let $R(X)$ be
the label of the chamber containing the corresponding part of the
upper arrangement, and let $S(X)$ be the label of the chamber
containing the corresponding part of the lower arrangement
(these corresponding parts can be obtained by overlaying the ansatz
arrangement with the upper and lower arrangements
respectively).
We label $X$ with the minor $\Delta^{R(X)}_{S(X)}$. The resulting
labeled ansatz arrangement for our example is given in
Figure~\ref{f:ansatzarrangement}.
\begin{figure}[htbp]
\beginpicture

\setcoordinatesystem units <0.35cm,0.35cm>             
\setplotarea x from -14 to 18, y from 0 to 8       
\linethickness=0.5pt           

\setlinear \plot 0 0  7 0  / %
\setlinear \plot 7 0  8 2  / %
\setlinear \plot 8 2  16 2  / %
\setlinear \plot 16 2 16.35 2.7  / %
\setlinear \plot 16.65 3.3 17 4  / %
\setlinear \plot 17 4  18 4  / %

\setlinear \plot 0 2  4 2  / %
\setlinear \plot 4 2  5 4  / %
\setlinear \plot 5 4  10 4 / %
\setlinear \plot 10 4 11 6  / %
\setlinear \plot 11 6 13 6  / %
\setlinear \plot 13 6 13.35 5.3  / %
\setlinear \plot 13.65 4.7  14 4  / %
\setlinear \plot 14 4  16 4  / %
\setlinear \plot 16 4  17 2  / %
\setlinear \plot 17 2  18 2  / %

\setlinear \plot 0 4  1 4  / %
\setlinear \plot 1 4  1.35 4.7  / %
\setlinear \plot 1.65 5.3  2 6  / %
\setlinear \plot 2 6  10 6  / %
\setlinear \plot 10 6 11 4  / %
\setlinear \plot 11 4 13 4  / %
\setlinear \plot 13 4 14 6  / %
\setlinear \plot 14 6 18 6  / %

\setlinear \plot 0 6  1 6  / %
\setlinear \plot 1 6  2 4  / %
\setlinear \plot 2 4  4 4  / %
\setlinear \plot 4 4  5 2  / %
\setlinear \plot 5 2  7 2  / %
\setlinear \plot 7 2  8 0  / %
\setlinear \plot 8 0  18 0  / %

\put{$\bullet$}[c] at 4.5 3
\put{$\bullet$}[c] at 7.5 1
\put{$\bullet$}[c] at 10.5 5

\put{$\Delta^{123}_{123}$}[c] at 0.3 5
\put{$\Delta^{124}_{124}$}[c] at 6 5
\put{$\Delta^{124}_{134}$}[c] at 12 5
\put{$\Delta^{123}_{134}$}[c] at 16 5
\put{$\Delta^{12}_{12}$}[c] at 2 3
\put{$\Delta^{12}_{14}$}[c] at 10.5 3
\put{$\Delta^{13}_{34}$}[c] at 17.5 3
\put{$\Delta^{1}_{1}$}[c] at 3.5 1
\put{$\Delta^{1}_{4}$}[c] at 13 1

\put{$1$}[c] at -0.7 0
\put{$2$}[c] at -0.7 2
\put{$3$}[c] at -0.7 4
\put{$4$}[c] at -0.7 6

\put{$t_2$}[c] at 4.5 -1
\put{$t_3$}[c] at 7.5 -1
\put{$m_4$}[c] at 10.5 -1

\endpicture
\caption{Ansatz arrangement for
$\underline{s_3} s_2 s_1 \underline{s_3 s_2}$.
Note that $g=\dot s_3 y_2(t_2) y_1(t_3) x_3(m_4) \dot s_3^{-1} \dot s_2$.}
\label{f:ansatzarrangement}
\end{figure}

Next, we note that the singular points in the ansatz arrangement
correspond precisely to the factors of $g$ of the form $x_{i_k}(m_k)$
and $y_{i_k}(t_k)$. We label these (beneath the arrangment) with
$t_k$ and $m_k$, respectively, for convenience.
The ansatz arrangement can then be used to compute the
coefficients $t_k$ and $m_k$ as follows.
Suppose $k\in J^-_{\mathbf{v}}\cup J^{\circ}_{\mathbf{v}}$.
Let $A_k$, $B_k$, $C_k$ and $D_k$ be the minors labelling the chambers
surrounding the singular point in the ansatz arrangement corresponding to
$k$, with $A_k$ and $D_k$ above and below it, and $B_k$ and $C_k$ on the
same horizontal level (see Figure~\ref{f:surroundingchambers}).
It is easy to check that Theorem~\ref{t:ansatz} implies that,
for $k\in J^{\circ}_{\mathbf{v}}$,
\begin{equation*}
t_k=\frac{A_k(z)D_k(z)}{B_k(z)C_k(z)},
\end{equation*}
and, for $k\in J^-_{\mathbf{v}}$,
\begin{equation*}
r_k=\frac{B_k(z)C_k(z)}{A_k(z)D_k(z)}.
\end{equation*}

\begin{figure}[htbp]
\beginpicture

\setcoordinatesystem units <0.35cm,0.35cm>             
\setplotarea x from -20 to 4, y from -1.5 to 4       
\linethickness=0.5pt           

\setlinear \plot 0 0  1 0 / %
\setlinear \plot 1 0  2 2 / %
\setlinear \plot 2 2  3 2 / %
\setlinear \plot 0 2  1 2 / %
\setlinear \plot 1 2  2 0 / %
\setlinear \plot 2 0  3 0 / %

\put{$\bullet$}[c] at 1.5 1

\put{$A_k$}[c] at 1.5 3
\put{$D_k$}[c] at 1.5 -1
\put{$B_k$}[c] at 0 1
\put{$C_k$}[c] at 3 1

\endpicture
\caption{Chambers surrounding the singular point corresponding to
$k\in J^-_{\mathbf{v}}\cup J^{\circ}_{\mathbf{v}}$.}
\label{f:surroundingchambers}
\end{figure}

\section{Determining Deodhar components for $SL_d$}
\label{s:componentexample}

\subsection{Graphical algorithm for Deodhar components}
In this section we show that Proposition~\ref{p:algorithm} also gives rise to
a graphical algorithm for $SL_d$ to determine the Deodhar component of
an element $B\in \mathcal{B}$.
Suppose that $B=z\dot w\cdot B^+$, where $w\in W$ and $z\in U^+$, and that
a reduced expression $\mathbf{w}$ for $w$ is chosen. Then
the graphical algorithm for computing the distinguished subexpression
$\mathbf{v}$ of $\mathbf{w}$ such that
$B\in {\mathcal R}_{\mathbf{v},\mathbf{w}}$ is as follows.

Firstly we draw the usual pseudoline diagram for the reduced
expression $\mathbf{w}$ of $w$ as in~\cite{BeFoZe:TotPos}
(call this the
{\em classical pseudoline arrangement} for $\mathbf{w}$). Each factor
$s_{i_k}$ of $\mathbf{w}$ corresponds to a singular crossing between the
$i_k$th and $i_{k+1}$st pseudolines from the bottom (see below for an
example). We define $v_{(0)}$ to be $1$.
Suppose that $v_{(0)}, v_{(1)}, \ldots ,v_{(k-1)}$ have already been
computed, and that if $k>1$
we have drawn the upper arrangement for the pair $\mathbf{v}_{(k-1)}=
(v_{(1)},v_{(2)},\ldots ,v_{(k-1)}),
\mathbf{w}_{(k-1)}=(w_{(0)},w_{(1)},\ldots ,w_{(k-1)})$.

We compute $v_{(k)}$ in the following way.
Consider the upper arrangement for the pair
$\mathbf{v}_{(k-1)},\mathbf{w}_{(k-1)}$.
If $k>1$, let $R_k$ be the label of the unbounded chamber at the
right hand end of this arrangement between the $i_k$th and $i_{k+1}$st
pseudolines (counting from bottom to top).
If $k=1$, we take $R_1$ to be $\{1,2,\ldots ,i_k\}$.

Let $S_k$ be the label of the chamber in the classical arrangement for
$\mathbf{w}$ immediately to the right of the $k$th crossing. It is
clear that $|R_k|=|S_k|=i_k$.
Now the minor $\Delta^{R_k}_{S_k}(z)$ determines the value of
$v_{(k)}$. We have, by Proposition~\ref{p:algorithm}:
\begin{itemize}
\item[$(a)$] If $v_{(k-1)}s_{i_k}>v_{(k-1)}$ and
$\Delta^{R_k}_{S_k}(z)\not=0$, then $v_{(k)}=v_{(k-1)}$.
\item[$(b)$] If $v_{(k-1)}s_{i_k}>v_{(k-1)}$ and
$\Delta^{R_k}_{S_k}(z)=0$, then $v_{(k)}=v_{(k-1)}s_{i_k}$.
\item[$(c)$] If $v_{(k-1)}s_{i_k}<v_{(k-1)}$, then
$v_{(k)}=v_{(k-1)}s_{i_k}$.
\end{itemize}
Thus $v_{(k)}$ is computed, and we draw the upper arrangement for the
pair $\mathbf{v}_{(k)},\mathbf{w}_{(k)}$, by building on the upper
arrangement for  $\mathbf{v}_{(k-1)},\mathbf{w}_{(k-1)}$ if $k>1$.
We are thus ready for the next step.

In this way, all of the $v_{(k)}$ are determined inductively. We also
note that at the end we have drawn the upper arrangement for the
pair $\mathbf{v},\mathbf{w}$. So, after drawing the lower arrangement for
$\mathbf{v},\mathbf{w}$, we are ready to apply the method in
Section~\ref{s:ansatztypea} to compute the factors of $g$ explicitly.

\subsection{An Example} \label{s:example}
\begin{figure}[htbp]
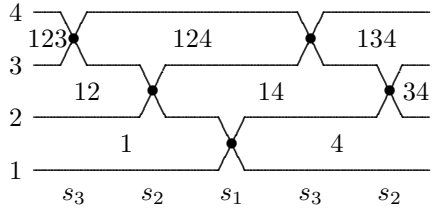

\beginpicture

\setcoordinatesystem units <0.35cm,0.35cm>             
\setplotarea x from -16 to 18, y from 0 to 8       
\linethickness=0.5pt           

\setlinear \plot 0 0  7 0  / %
\setlinear \plot 7 0  8 2  / %
\setlinear \plot 8 2  13 2  / %
\setlinear \plot 13 2 14 4  / %
\setlinear \plot 14 4 15 4  / %

\setlinear \plot 0 2  4 2  / %
\setlinear \plot 4 2  5 4  / %
\setlinear \plot 5 4  10 4 / %
\setlinear \plot 10 4 11 6  / %
\setlinear \plot 11 6 15 6  / %

\setlinear \plot 0 4  1 4  / %
\setlinear \plot 1 4  2 6  / %
\setlinear \plot 2 6  10 6  / %
\setlinear \plot 10 6 11 4  / %
\setlinear \plot 11 4 13 4  / %
\setlinear \plot 13 4 14 2  / %
\setlinear \plot 14 2 15 2  / %

\setlinear \plot 0 6  1 6  / %
\setlinear \plot 1 6  2 4  / %
\setlinear \plot 2 4  4 4  / %
\setlinear \plot 4 4  5 2  / %
\setlinear \plot 5 2  7 2  / %
\setlinear \plot 7 2  8 0  / %
\setlinear \plot 8 0  15 0  / %

\put{$\bullet$}[c] at 1.5 5
\put{$\bullet$}[c] at 4.5 3
\put{$\bullet$}[c] at 7.5 1
\put{$\bullet$}[c] at 10.5 5
\put{$\bullet$}[c] at 13.5 3

\put{$123$}[c] at 0.5 5
\put{$124$}[c] at 6 5
\put{$134$}[c] at 13 5
\put{$12$}[c] at 2 3
\put{$14$}[c] at 9 3
\put{$34$}[c] at 14.5 3
\put{$1$}[c] at 3.5 1
\put{$4$}[c] at 11.5 1

\put{$1$}[c] at -0.7 0
\put{$2$}[c] at -0.7 2
\put{$3$}[c] at -0.7 4
\put{$4$}[c] at -0.7 6

\put{$s_3$}[c] at 1.5 -1
\put{$s_2$}[c] at 4.5 -1
\put{$s_1$}[c] at 7.5 -1
\put{$s_3$}[c] at 10.5 -1
\put{$s_2$}[c] at 13.5 -1

\endpicture
\caption{Classical arrangement for $\mathbf{w}=s_3s_2s_1s_3s_2$.}
\label{f:classicalwarrangement}
\end{figure}

\begin{figure}[htbp]
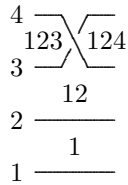

\beginpicture

\setcoordinatesystem units <0.35cm,0.35cm>             
\setplotarea x from -22 to 18, y from 0 to 8       
\linethickness=0.5pt           

\setlinear \plot  0 0  3 0 / %

\setlinear \plot  0 2  3 2 / %

\setlinear \plot 0 4  1 4  / %
\setlinear \plot 1 4  1.35 4.7  / %
\setlinear \plot 1.65 5.3  2 6 / %
\setlinear \plot 2 6  3 6  / %

\setlinear \plot 0 6  1 6 / %
\setlinear \plot 1 6  2 4 / %
\setlinear \plot 2 4  3 4 / %

\put{$123$}[c] at 0.3 5
\put{$124$}[c] at 2.7 5
\put{$12$}[c] at  1.5 3
\put{$1$}[c] at 1.5 1

\put{$1$}[c] at -0.7 0
\put{$2$}[c] at -0.7 2
\put{$3$}[c] at -0.7 4
\put{$4$}[c] at -0.7 6

\endpicture
\caption{Upper arrangement for
$\underline{s_3}$.
Note that $g_{(1)}=\dot s_3$.}
\label{f:arrangement1}
\end{figure}

\begin{figure}[htbp]
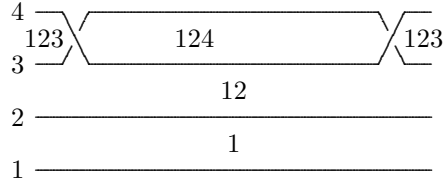

\beginpicture

\setcoordinatesystem units <0.35cm,0.35cm>             
\setplotarea x from -16 to 18, y from 0 to 8       
\linethickness=0.5pt           

\setlinear \plot  0 0  15 0 / %

\setlinear \plot  0 2  15 2 / %

\setlinear \plot 0 4  1 4  / %
\setlinear \plot 1 4  1.35 4.7  / %
\setlinear \plot 1.65 5.3  2 6 / %
\setlinear \plot 2 6  13 6  / %
\setlinear \plot 13 6 13.35 5.3  / %
\setlinear \plot 13.65 4.7  14 4  / %
\setlinear \plot 14 4  15 4  / %

\setlinear \plot 0 6  1 6 / %
\setlinear \plot 1 6  2 4 / %
\setlinear \plot 2 4  13 4 / %
\setlinear \plot 13 4 14 6 / %
\setlinear \plot 14 6 15 6 / %

\put{$123$}[c] at 0.3 5
\put{$124$}[c] at 6 5
\put{$123$}[c] at 14.7  5
\put{$12$}[c] at 7.5 3
\put{$1$}[c] at 7.5 1

\put{$1$}[c] at -0.7 0
\put{$2$}[c] at -0.7 2
\put{$3$}[c] at -0.7 4
\put{$4$}[c] at -0.7 6

\endpicture
\caption{Upper arrangement for
$\underline{s_3} s_2 s_1 \underline{s_3}$.
Note that $g_{(4)}=\dot s_3 y_2(t_2) y_1(t_3) x_3(m_4) \dot s_3^{-1}$.}
\label{f:arrangement4}
\end{figure}

In this section we give an explicit example of the graphical algorithm
described above. We consider the element
$$z=\left( \begin{array}{cccc} 1 & 1 & 2 & 1 \\ 0 & 1 & 4 & 2 \\
0 & 0 & 1 & 0 \\ 0 & 0 & 0 & 1 \end{array} \right)\in SL_4(\mathbb{C}),$$
and set $w=s_3s_2s_1s_3s_2$, so we have the reduced expression
$\mathbf{w}=(1,s_3,s_3s_2,s_3s_2s_1,s_3s_2s_1s_3,s_3s_2s_1s_3s_2)$ for $w$.
The classical arrangement for $\mathbf{w}$ is given in
Figure~\ref{f:classicalwarrangement}.
We start with $v_{(0)}=1$.
\begin{description}
\item[Step 1.] We have $v_{(0)}s_3=s_3>v_{(0)}$.
Here, $R_1=\{1,2,3\}$ and $S_1=\{1,2,4\}$.
Since $i_1=3$ and $\Delta^{123}_{124}(z)=0$,
we are in case (b) and $v_{(1)}=v_{(0)}s_3=s_3$.
\item[Step 2.] We have $v_{(1)}s_2=s_3s_2>v_{(1)}$.
The upper arrangement for the pair $\mathbf{v}_{(1)},\mathbf{w}_{(1)}$
is shown in Figure~\ref{f:arrangement1}.
Here, $R_2=\{1,2\}$ and $S_2=\{1,4\}$.
Since $i_2=2$ and $\Delta^{12}_{14}(z)=2\not=0$,
we are in case (a) and $v_{(2)}=v_{(1)}=s_3$.
\item[Step 3.] We have $v_{(2)}s_1=s_3s_1>v_{(2)}$.
The upper arrangement for the pair $\mathbf{v}_{(2)},\mathbf{w}_{(2)}$
is shown in Figure~\ref{f:arrangement2}.
Here, $R_3=\{1\}$ and $S_3=\{4\}$.
Since $i_3=1$ and $\Delta^{1}_{4}(z)=1$,
we are in case (a) and $v_{(3)}=v_{(2)}=s_3$.
\item[Step 4.] We have $v_{(3)}s_3=s_3s_3=1<v_{(2)}$.
The upper arrangement for the pair $\mathbf{v}_{(3)},\mathbf{w}_{(3)}$
is shown in Figure~\ref{f:arrangement3}.
We are in case (c) and $v_{(4)}=v_{(3)}s_3=1$.
\item[Step 5.] We have $v_{(4)}s_2=s_2>v_{(4)}$.
The upper arrangement for the pair $\mathbf{v}_{(4)},\mathbf{w}_{(4)}$
is shown in Figure~\ref{f:arrangement4}.
Here, $R_5=\{1,2\}$ and $S_5=\{3,4\}$.
Since $i_5=2$ and $\Delta^{12}_{34}(z)=0$,
we are in case (b) and $v_{(5)}=v_{(4)}s_2=s_2$.
The upper arrangement for the pair $\mathbf{v}_{(5)},\mathbf{w}_{(5)}$
is shown in Figure~\ref{f:exampleupperarrangement}.
\end{description}

\begin{figure}[htbp]
\beginpicture

\setcoordinatesystem units <0.35cm,0.35cm>             
\setplotarea x from -21 to 18, y from 0 to 8       
\linethickness=0.5pt           

\setlinear \plot  0 0  6 0 / %

\setlinear \plot  0 2  6 2 / %

\setlinear \plot 0 4  1 4  / %
\setlinear \plot 1 4  1.35 4.7  / %
\setlinear \plot 1.65 5.3  2 6 / %
\setlinear \plot 2 6  6 6  / %

\setlinear \plot 0 6  1 6 / %
\setlinear \plot 1 6  2 4 / %
\setlinear \plot 2 4  6 4 / %

\put{$123$}[c] at 0.3 5
\put{$124$}[c] at 4 5
\put{$12$}[c] at  3 3
\put{$1$}[c] at 3 1

\put{$1$}[c] at -0.7 0
\put{$2$}[c] at -0.7 2
\put{$3$}[c] at -0.7 4
\put{$4$}[c] at -0.7 6

\endpicture
\caption{Upper arrangement for
$\underline{s_3} s_2 $.
Note that $g_{(2)}=\dot s_3 y_2(t_2)$.}
\label{f:arrangement2}
\end{figure}

\begin{figure}[htbp]
\beginpicture

\setcoordinatesystem units <0.35cm,0.35cm>             
\setplotarea x from -19.5 to 18, y from 0 to 8       
\linethickness=0.5pt           

\setlinear \plot  0 0  9 0 / %

\setlinear \plot  0 2  9 2 / %

\setlinear \plot 0 4  1 4  / %
\setlinear \plot 1 4  1.35 4.7  / %
\setlinear \plot 1.65 5.3  2 6 / %
\setlinear \plot 2 6  9 6  / %

\setlinear \plot 0 6  1 6 / %
\setlinear \plot 1 6  2 4 / %
\setlinear \plot 2 4  9 4 / %

\put{$123$}[c] at 0.3 5
\put{$124$}[c] at 5.5 5
\put{$12$}[c] at  4.5 3
\put{$1$}[c] at 4.5 1

\put{$1$}[c] at -0.7 0
\put{$2$}[c] at -0.7 2
\put{$3$}[c] at -0.7 4
\put{$4$}[c] at -0.7 6

\endpicture
\caption{Upper arrangement for
$\underline{s_3} s_2 s_1$.
Note that $g_{(3)}=\dot s_3 y_2(t_2) y_1(t_3)$.}
\label{f:arrangement3}
\end{figure}

We conclude that $B$ lies in the Deodhar component
${\mathcal R}_{\mathbf{v},\mathbf{w}}$, where
$\mathbf{v}=(1,s_3,s_3,s_3,1,s_2)$ and
$\mathbf{w}=(1,s_3,s_3s_2,s_3s_2s_1,s_3s_2s_1s_3,s_3s_2s_1s_3s_2)$.

We remark that we can use the above computation to determine criteria for
$z\dot w\cdot B^+$ to lie in $\mathcal{R}_{\mathbf{v},\mathbf{w}}$,
where $z$ is an arbitrary matrix in $U^+$, in terms of minors of $z$
(see Corollary~\ref{c:DeoInequalities}). Suppose that
$$z=\left( \begin{array}{cccc} 1 & a_{12} & a_{13} & a_{14} \\
0 & 1 & a_{23} & a_{24} \\
0 & 0 & 1 & a_{34} \\ 0 & 0 & 0 & 1 \end{array} \right).$$
We obtain that
\begin{equation*}
\mathcal{R}_{\mathbf{v},\mathbf{w}}=
\{z\dot w\cdot B^+\,:\,
\Delta^{123}_{124}(z)=a_{34}=0,
\Delta^{12}_{14}(z)=a_{24}\not=0, \Delta^{1}_{4}(z)=a_{14}\not=0,
\Delta^{12}_{34}(z)=a_{13}a_{24}-a_{14}a_{23}=0
\}.
\end{equation*}
We note that no condition is obtained on the entry $a_{12}$ of $z$.
No such condition could arise, since $z\dot w\cdot B^+$ doesn't depend on
$a_{12}$, as $w\inv\alpha_1$ is positive.

Finally, we note that in our example, $B$ lies in the Deodhar
component $\mathcal{R}_{\mathbf{v},\mathbf{w}}$ considered in
Section~\ref{s:ansatztypea}, so we can apply the graphical algorithm given
there in order to compute the coefficients $t_k$ and $m_k$ in the factorisation
of $g=\dot s_3 y_2(t_2) y_1(t_3) x_3(m_4) \dot s_3^{-1} \dot s_2\in
G_{\mathbf{v},\mathbf{w}}$ (where $B=g\cdot B^+$). We obtain:
\begin{eqnarray*}
t_2 & = & \frac{\Delta^{124}_{124}(z)\Delta^{1}_{1}(z)}
{\Delta^{12}_{12}(z)\Delta^{12}_{14}(z)}=1/2, \\
t_3 & = & \frac{\Delta^{12}_{14}(z)\Delta^{\emptyset}_{\emptyset}(z)}
{\Delta^{1}_{1}(z)\Delta^{1}_{4}(z)}=2,
\end{eqnarray*}
and
$$
r_4 = \frac{\Delta^{124}_{124}(z)\Delta^{124}_{134}(z)}
{\Delta^{1234}_{1234}(z)\Delta^{12}_{14}(z)}=2. \\
$$

Finally, we show that the minor
$\Delta^{v_{(k-1)}\omega_{i_k}}_{s_{i_k}\omega_{i_k}}(g_{(k-1)})$,
that would appear as correction term vanishes on this Deodhar component,
so that $m_k=r_k$. We have
$\Delta^{v_{(k-1)}\omega_{i_k}}_{s_{i_k}\omega_{i_k}}(g_{(k-1)})=
\Delta^{v_{(3)}\omega_{i_4}}_{s_{i_4}\omega_{i_4}}(g_{(3)})=
\Delta^{s_3\omega_3}_{s_3\omega_3}(\dot s_3 y_2(t_2)y_1(t_3))$.
We note that $g_{(3)}\dot s_3\omega_3$ is a linear combination of elements
of weight $\omega_3$ and
$s_3(s_3\omega_3-\alpha_2)=\omega_3-\alpha_2-\alpha_3$,
and therefore has zero component in the $s_3\omega_3$-weight space, from
which it follows that
$\Delta^{s_3\omega_3}_{s_3\omega_3}(\dot s_3 y_2(t_2)y_1(t_3))=0$,
so $m_4=r_4-0=r_4$. Thus, in this case, $m_4=2$.
We finally see that
$$\left( \begin{array}{cccc} 1 & 1 & 2 & 1 \\ 0 & 1 & 4 & 2 \\
0 & 0 & 1 & 0 \\ 0 & 0 & 0 & 1 \end{array} \right)\dot s_3\dot s_2\dot s_1
\dot s_3\dot s_2\cdot B^+
=
\dot s_3 y_2(1/2) y_1(2) x_3(2) \dot s_3^{-1} \dot s_2\cdot B^+.$$
\begin{rem}
It is possible that the minor
$\Delta^{v_{(k-1)}\omega_{i_k}}_{s_{i_k}\omega_{i_k}}(g_{(k-1)})$
is non-zero. This happens, for example, for $G=SL_4$. Let $v=1$ and
$w=w_0$, and take
$\mathbf{v}=(1,1,1,1,s_1,s_1,1)$, a subexpression for $v$ in
the reduced expression for $w$ with factors $(s_1,s_2,s_3,s_1,s_2,s_1)$.
Then $g_{(5)}=y_1(t_1)y_2(t_2)y_3(t_3)\dot s_1 y_2(t_5)$, and
it is easy to check that
$$\Delta^{v_{(5)}\omega_{i_6}}_{s_{i_6}\omega_{i_6}}(g_{(5)})=-t_1.$$
\end{rem}
\section{Total positivity}\label{s:totpos}
>From now on let $\k=\R$. We view the group $G$ and the flag
variety $\mathcal B$ as real manifolds, with the corresponding
Hausdorff topology.

\begin{defn}[\cite{Lus:TotPos94}]
The totally nonnegative part $U^-_{\ge 0}$ of $U^-$ is defined
to be the semigroup in $U^-$ generated by the $y_i(t)$ for $t\in
\R_{\ge 0}$. The totally nonnegative part of $\mathcal B$ is
defined by
 $$
\mathcal B_{\ge 0}:=\overline{\{u\cdot B^+\ |\ u\in U^-_{\ge 0}\}}
 $$
where the closure is taken inside $\mathcal B$ in its real
topology.
\end{defn}

Let us collect below some useful facts, see
\cite{Lus:TotPos94} and also \cite{BeZel:TotPos} for $(1)$.
\begin{enumerate}
\item
For any braid relation such as $s_i s_j s_i=s_j s_i s_j$ in $W$
there is a subtraction-free rational transformation relating the
parameters of the corresponding parametrizations.
For example in type $A_2$,
 $$
y_i(a)y_j(b)y_i(c)= y_j\left(\frac{bc}{a+c}\right) y_i(a+c)\,
y_j\left(\frac{ab}{a+c}\right).
 $$
\item For $w\in W$ and a reduced expression $w=s_{i_1}\dotsc s_{i_n}$ define
 $$
U^-_{> 0}(w):=\{y_{i_1}(t_1)\cdot \dotsc \cdot y_{i_n}(t_n) \ |\
t_1,\dotsc,t_n\in \R_{>0}\}.
 $$
By $(1)$ this set is independent of the reduced expression chosen.
Moreover any product of $y_i(t)$'s for positive parameters $t$ can
be transformed until it is seen to lie in some
$U^-_{>0}(w)$. Therefore
 $$
 U^-_{\ge 0}=\bigsqcup_{w\in W}U^-_{> 0}(w).
 $$
This is of course precisely the decomposition of $ U^-_{\ge
0}$ induced by Bruhat decomposition, that is,  $
U^-_{>0}(w)=U^-_{\ge 0}\cap B^+\dot w B^+$.
\item
The totally positive part of $U^-$ may be defined as
$U^-_{>0}=U^-_{\ge 0}(w_0)$. For the flag variety the totally
positive part is taken to be
 $$\mathcal B_{>0}:=\{u\cdot B^+\ |\ u\in U^-_{>0}\}.$$
Clearly $\mathcal B_{>0}$ is open dense in $\mathcal B_{\ge 0}$.
\item Let $u\in U^-_{\ge 0}$. The semigroup property $uU^-_{\ge
0}\subset U^-_{\ge 0}$ implies, by continuity, that
 $$u\cdot\mathcal B_{\ge 0}\subset \mathcal B_{\ge 0}.$$
\end{enumerate}

\begin{defn} For $v,w\in W$ with $v\le w$, let
 $$
\mathcal R_{v,w}^{>0} :=\mathcal R_{v,w}\cap \mathcal B_{\ge 0}.
 $$
\end{defn}

In the special case where $v=1$ we have $\mathcal R_{1,w}\cong
U^-\cap B^+\dot w B^+$ and $\mathcal R_{1,w}^{>0}=U^-_{>
0}(w)\cdot B^+$ (see property $(2)$ above). From this observation
Lusztig \cite{Lus:TotPos94} conjectured that also $\mathcal
R_{v,w}^{>0}$ is a semi-algebraic cell. The first proof of this is
in \cite{Rie:CelDec}. However this proof does not provide an
explicit parametrization and uses deep properties of canonical
bases. We will now give a different proof which gives
parametrizations and is completely elementary.

Let us choose a reduced expression $\w$ for $w$ with factors
$(s_{i_1},\dotsc,s_{i_n})$. To $v\le w$ we may associate the
positive subexpression $\v_{+}$ for $v$ in $\w$ as in
Lemma~\ref{l:positive}. Note that $\v_{+}$ is non-decreasing, so
$J_{\v_{+}}^-=\emptyset$. We define
 $$
G_{\v_{+},\w}^{>0}:=\left\{g=g_1 g_2\cdots g_n\,\left |\,
\begin{array}{ll} g_k=
y_{i_k}(t_k)
\text{ for $t_k\in\R_{>0}$,}&\text{ if $k\in J^\circ_{\v_{+}}$}\\
g_k=\dot s_{i_k}, &\text{ if $k\in J^+_{\v_{+}}$}
\end{array} \right .
\right\}
 $$
Then $G_{\v_{+},\w}^{>0}\cong \R_{>0}^{\ell(w)-\ell(v)}$ is a
semi-algebraic cell in $G$.  The aim of this section is to prove
the following theorem.
\begin{thm}\label{t:totpos}
The isomorphism $G_{\v_{+},\w}\overset\sim\to\mathcal
R_{\v_{+},\w}$ restricts to an isomorphism of real semi-algebraic
varieties
 $$
G^{>0}_{\v_{+},\w}\overset\sim\To\mathcal R_{v,w}^{>0}.
 $$
\end{thm}

Note that if $v=1$ then $G^{>0}_{\v_+,\w}=U^-_{>0}(w)$ and, as a subset
of $G$, does not depend on the reduced expression $\w$. This is no
longer true if $v\ne 1$ as can be seen already in type $A_2$. We
begin the proof of Theorem~\ref{t:totpos} with a simple
observation about minors.
\begin{lem}\label{l:posminors}
Suppose $B=z\dot w\cdot B^+$ lies in $\mathcal B_{\ge 0}$ with
$z\in U^+$ and $w\in W$. For any dominant weight $\lambda$ and
$v\in W$,
 $$
 \Delta_{w\lambda}^{v\lambda}(z)\ge 0.
 $$
\end{lem}
\begin{proof} Since $B\in\mathcal B_{\ge 0}$ we can find a sequence $u_n\cdot
B^+$ with $u_n\in U^-_{>0}$ that converges to $B=z\dot w\cdot
B^+$. Note that for any $u=y_{i_1}(t_1)\dotsc y_{i_N}(t_N)\in U^-_{>0}$
and $x\in W$, the element $u\cdot \xi_\lambda$ in $V(\lambda)$
has a positive projection to the $x\lambda$ weight space, using that
$s_{i_1}\dotsc s_{i_N}=w_0$ has a subexpression for $x$. Now
we have
 $$
\frac{u_n\cdot \xi_\lambda}{\left<u_n\cdot \xi_\lambda, \dot
w\cdot\xi_\lambda\right>}\to z\dot w\cdot\xi_\lambda \qquad\qquad
(n\to \infty),
 $$
where the denominator $\left<u_n\cdot \xi_\lambda, \dot
 w\cdot\xi_\lambda\right>$
is just the required normalization factor.
It follows that
 $$
0\le \lim_{n\to\infty}\frac{ \left<u_n\cdot \xi_\lambda, \dot
 v\cdot\xi_\lambda\right>}{\left<u_n\cdot \xi_\lambda, \dot
 w\cdot\xi_\lambda\right> }=\left<z\dot
w\cdot\xi_\lambda,\dot v\cdot\xi_\lambda \right> =
\Delta_{w\lambda}^{v\lambda}(z).
 $$

\end{proof}

We need to recall one lemma.

\begin{lem}[\cite{Rie:CelDec}~Lemma~2.3]
\label{l:oldlem} Suppose $w=w_1 w_2$ with
$\ell(w)=\ell(w_1)+\ell(w_2)$. Consider the reduction map
$\pi^w_{w_1}: B^+\dot w\cdot B^+\to B^+\dot w_1\cdot B^+ $. If
$B\in B^+\dot w\cdot B^+$ lies in $\mathcal B_{\ge 0}$ then so
does $\pi^w_{w_1}(B)$.
\end{lem}

This lemma is easy to see if $B\in \mathcal R_{1,w}^{>0}$~: In
that case using Lusztig's parametrization we may write
$B=y_{i_1}(t_1)\dotsc y_{i_n}(t_n)\cdot B^+$ for some positive
$t_i$, being careful to choose a reduced expression $s_{i_1}\cdots
s_{i_n}$ of $w$ such that $w_1=s_{i_1}\cdots s_{i_m}$, where
$m=\ell(w_1)$. The element $\pi^w_{w_1}(B)=y_{i_1}(t_1)\dotsc
y_{i_m}(t_m)\cdot B^+$ is then clearly totally nonnegative again.
The property extends from the dense open part $\mathcal R_{1,w}$
to the whole Bruhat cell essentially by continuity (see
\cite{Rie:CelDec} for a more careful argument).

Now we are ready to show one part of the theorem.
\begin{lem}\label{l:TotPosThenPos}
If $\mathcal R_{\v,\w}\cap \mathcal B_{\ge 0}\ne\emptyset$, then
$\v$ is a positive subexpression of $\w$.
\end{lem}
\begin{proof}
Let $B\in\mathcal R_{\v,\w}\cap\mathcal B_{\ge 0}$ and write
$B=z\dot w\cdot B^+$ for $z\in U^+$. Suppose $\v\prec\w$ is not a
positive subexpression. Then $J_\v^{-}\ne\emptyset$. Let $k\in
J_\v^{-}$. The equation \eqref{e:generalansatz} together with
Remark~\ref{r:formulas}.(1) gives
 $$
-1=\frac{
\prod_{j\not=i_k}\Delta_{w_{(k)}\omega_j}^{v_{(k)}\omega_j}(z)^{-a_{j,i_k}}}
{\Delta_{w_{(k)}\omega_{i_k}}^{v_{(k)}\omega_{i_k}}(z)
\Delta_{w_{(k-1)}\omega_{i_k}}^{v_{(k-1)}\omega_{i_k}}(z)}.
 $$
Therefore at least one of the minors in this formula must be negative. By
Lemma~\ref{l:posminors} this implies that one of the two elements
$z\dot w_{(k)}\cdot B^+$ and $z\dot w_{(k-1)}\cdot B^+$ does not
lie in $\mathcal B_{\ge 0}$. Since these are both reductions of
$B$ we have a contradiction to Lemma~\ref{l:oldlem}.
\end{proof}

\begin{rem}\label{r:posminors}
Suppose $z\dot w\cdot B^+\in \mathcal R_{v,w}^{>0}$ and
$\w$ is a reduced expression for $w$ with positive subexpression $\v_+$
for $v$. Then by a combination of the above lemmas we have,
 $$
\Delta^{v_{(k)}\omega_{i_k}}_{w_{(k)}\omega_{i_k}}(z)>0,\qquad \quad
 k=0,1,\dotsc,\ell(w).
$$
Recall that these minors, as standard chamber minors, are
automatically nonzero, hence the strict positivity.
Since $\v_+$ is non-decreasing
there are no special chamber minors.
So the observation is that if $z\dot w\cdot B^+$
lies in $\mathcal B_{\ge 0}$, then all of the
associated chamber minors are positive.
\end{rem}

The following lemma is a technical tool we will need to finish the
proof of the theorem.

\begin{lem}\label{l:tool}
Let $v\le w$ in $W$ and suppose $\alpha_{i_0}$ is a simple root
such that
$u\inv\alpha_{i_0}>0$ for all $v\le u\le w$.  Then for all $g\cdot
B^+\in \mathcal R_{v,w}$ and any $m\in \R$,
\begin{equation}\label{e:invariance}
  x_{i_0}(m)g\cdot B^+=g\cdot B^+.
\end{equation}
In other words, if $u\inv\alpha_{i_0}>0$ for all $v\le u\le w$, then
$\mathcal R_{v,w}$ is contained in the Springer fiber of $x_{i_0}(m)$.
\end{lem}
Note that it is easy to see that the condition on $\alpha_{i_0}$
is also necessary. Suppose $m\ne 0$. If $x_{i_0}(m)$ fixes the
elements of $\mathcal R_{v,w}$ then it also fixes the elements of
the closure. So in particular, $x_{i_0}(m)\dot u\cdot B^+=\dot
u\cdot B^+$ for $v\le u\le w$. This implies $u\inv\alpha_{i_0}>0$.
\begin{proof}
Let $\v_+$ be the positive subexpression for $v$ of a reduced
expression $\w$ for $w$. Since the corresponding Deodhar component
$\mathcal R_{\v_+,\w}$ is dense in $\mathcal R_{v,w}$, it suffices
to show that $x_{i_0}(m)g\cdot B^+=g\cdot B^+$ for $g\cdot
B^+\in\mathcal R_{\v_+,\w}$. In other words we may assume $g\in
G_{\v_+,\w}$.

By the defining property
\eqref{e:PositiveSubexpression} for positive subexpressions
we have that
 $
 v_{(j-1)}\alpha_{i_j}>0
 $
for all $j$. Also $J^-_{\v_+}=\emptyset$. So we may write $g\in G_{\v_+,\w}$ as
 $$
 g=\left(\prod_{j\in J^{\circ}_{\v_+}}
 y_{v_{(j-1)}\alpha_{i_j}}(t_j)\right) \dot v,
 $$
where $y_{v_{(j-1)}\alpha_{i_j}}(t):=\dot v_{(j-1)}y_{i_j}(t)\dot
v_{(j-1)}\inv\in U^-$. Let us set
$y=\prod_{j\in J^{\circ}_{\v_+}} y_{v_{(j-1)}\alpha_{i_j}}(t_j)$. Then
we have
 $$
 g\cdot B^+=y \dot v\cdot B^+=z\dot w\cdot B^+,
 $$
where $z\in U^+$. Let us also write the reductions
$\pi^w_{w_{(k)}}(g\cdot B^+)$ as
 $$
 g_{(k)}\cdot B^+=y_{(k)} \dot v_{(k)}\cdot B^+=z\dot w_{(k)}\cdot B^+,
 $$
where $y_{(k)}=\prod_{j\in J^{\circ}_{\v_+}, j\le
k } y_{v_{(j-1)}\alpha_{i_j}}(t_j)$ and otherwise the notation is
as usual.

\vskip .1cm
We will now show that the conditions on $\alpha_{i_0}$ imply the following
assertions.

\begin{enumerate}
\item[$(i)$]
If $x\in W$ satisfies $v_{(j)}\le x\le w_{(j)}$ for some $j=1,\dotsc, n$,
then $x\inv \alpha_{i_0}>0$.
\item[$(ii)$]
$
v_{(j-1)}\alpha_{i_j}\ne \alpha_{i_0}
$
for all $j=1,\dotsc, n$.
\end{enumerate}

\vskip .1cm
For $x\in W$ let $R^+(x):=\{\alpha\in R^+\, |\, x\inv \alpha<0\}$.
Suppose we have an $x\in W$ with $v_{(j)}\le x\le w_{(j)}$.
We want to show that there exists $x'$ with $v\le x'\le w$
such that $R^+(x)\subseteq R^+(x')$. This will imply $(i)$.
Set $x_{(j)}:=x$. It suffices if we can construct $x_{(j+1)}$ with
$v_{(j+1)}\le x_{(j+1)}\le w_{(j+1)}$ and
$R^+(x_{(j)})\subseteq R^+(x_{(j+1)})$.
For this there are two cases.
If already $v_{(j+1)}\le x_{(j)}$, then we may set
$x_{(j+1)}:=x_{(j)}$. Otherwise we
must be in the situation $v_{(j+1)}=v_{(j)}s_{i_{j+1}}$ and we need to take
(at least) $x_{(j+1)}:=x_{(j)}s_{i_{j+1}}$ to obtain
$v_{(j+1)}\le x_{(j+1)}\le w_{(j+1)}$.
Because $v_{(j+1)}\le x_{(j+1)}$ and $v_{(j+1)}\not \le x_{(j)}$, we find that
$x_{(j)}<x_{(j+1)}$ and also
$R^+(x_{(j)})\subset R^+(x_{(j+1)})$. So $x_{(j+1)}$ has been
constructed successfully.

Now consider $x=v_{(j-1)}s_{i_j}$.
Clearly $v_{(j)}\le x\le w_{(j)}$ is satisfied, and so by $(i)$ we have
$s_{i_j}\inv v_{(j-1)}\inv\alpha_{i_0}>0$. Therefore
 $$
\ell(s_{i_0}v_{(j-1)}s_{i_j})=
\ell(v_{(j-1)}s_{i_j})+1=\ell(v_{(j-1)})+2
 $$
using also that $\v_+$ is non-decreasing.
It follows that $s_{i_0}v_{(j-1)}s_{i_j}>s_{i_0}v_{(j-1)}$ and
$ s_{i_0}v_{(j-1)}\alpha_{i_j}>0 $.
This implies $(ii)$.

\vskip .1cm

Finally we can put everything together to
show that $x_{i_0}(m)g\cdot B^+$ lies in $\mathcal
R_{\v_+,\w}$ and equals $g\cdot B^+$. By
Corollary~\ref{c:DeoInequalities} and Theorem~\ref{t:ansatz} we
know exactly which minors to check. Namely we only have to prove
\begin{enumerate}
\item
 $$
 \Delta^{v_{(k-1)}\omega_{i_k}}_{w_{(k)}\omega_{i_k}}(x_{i_0}(m)z)=0,
 \qquad  k\in J^+_{\v_+},
 $$
\item
 $$
 \Delta^{v_{(k)}\omega_{i_k}}_{w_{(k)}\omega_{i_k}}(x_{i_0}(m)z)=
 \Delta^{v_{(k)}\omega_{i_k}}_{w_{(k)}\omega_{i_k}}(z), \qquad
  k=1,\dotsc, n.
 $$
\end{enumerate}
Suppose $l\in\Z_{\ge 0}$ and
$\mu={v_{(k-1)}\omega_{i_k}-l\alpha_{i_0}}$. Let
$\zeta=pr_{\mu}(y_{(k)}\dot v_{(k)}\cdot\xi_{\omega_{i_k}})$ in
$V(\omega_{i_k})$. If
$\zeta\ne 0$ then the weight $\mu$ must be of the form
 $$
 v_{(k-1)}\omega_{i_k}-l\alpha_{i_0}=v_{(k)}\omega_{i_k}-
 \sum_{j\le k, j\in J^\circ_{\v}} c_j v_{(j-1)}\alpha_{i_j}
 $$
with $c_j\ge 0$, since the factors of $y_{(k)}$ are
$y_{v_{(j-1)}\alpha_{i_j}}(t_j)$ for $j\in J^{\circ}_{\v_+}$ with
$j\le k$.
Simplifying this equation we get
 $$
l\alpha_{i_0}=v_{(k-1)}\alpha_{i_k}+
 \sum_{j\le k, j\in J^\circ_{\v_+}} c_j v_{(j-1)}\alpha_{i_j}.
 $$
But the right hand side is a non-zero sum of positive roots
$\alpha$ not equal to $\alpha_{i_0}$ by $(ii)$ above.
Therefore we have a contradiction and so $\zeta=0$. Since $z\dot
w_{(k)}\cdot\xi_{\omega_{i_k}}$ and $y_{(k)}\dot v_{(k)}\cdot
\xi_{\omega_{i_k}}$ are collinear in $V(\omega_{i_k})$ this
implies that also
 $$
 pr_{v_{(k-1)}\omega_{i_k}-l\alpha_{i_0}}
(z \dot w_{(k)}\cdot\xi_{\omega_{i_k}})=0.
 $$
for all $l\ge 0$. Therefore $\left< x_{i_0}(m)z\dot
w_{(k)}\cdot\xi_{\omega_{i_k}}\ ,\ \dot v_{(k-1)}\cdot
\xi_{\omega_{i_k}} \right>=0$ and $(1)$ holds.

We can now prove $(2)$ in a completely analogous way. Let
$l\in\Z_{>0}$ and consider
$\zeta=pr_{v_{(k)}\omega_{i_k}-l\alpha_{i_0}}(y_{(k)}\dot
v_{(k)}\cdot\xi_{\omega_{i_k}})$. Then $\zeta\ne 0$ only if
 $$
 l\alpha_{i_0}=\sum_{j\le k, j\in J^\circ_{\v_+}} c_j v_{(j-1)}\alpha_{i_j}
 $$
for $c_j\ge 0$. And again this is impossible by
$(ii)$. So $\zeta=0$ and with it
 $$
 pr_{v_{(k)}\omega_{i_k}-l\alpha_{i_0}}(z\dot w_{(k)}\cdot\xi_{\omega_{i_k}})=0
 $$
for all $l>0$. This implies $(2)$.
\end{proof}

\begin{lem}\label{l:welldef}
If $g\in G^{>0}_{\v_+,\w}$ then $g\cdot B^+\in \mathcal
R_{v,w}^{>0}$.
\end{lem}

\begin{proof}
By Proposition~\ref{p:parameterization} we have $g\cdot B^+\in
\mathcal R_{\v_+,\w}$. We need to show that $g\cdot B^+\in\mathcal
B_{\ge 0}$. We have $g=g_1\dotsc g_{n}$ with
 $$
g_j=\begin{cases}
y_{i_j}(t_j) \text{ for $t_j\in\R_{>0}$,}&
\text{ if $j\in J^\circ_{\v_{+}}$}\\
\dot s_{i_j}, &\text{ if $j\in J^+_{\v_{+}}$}\end{cases}
 $$
Clearly $g_{n}\cdot B^+\in \mathcal B_{\ge 0}$. We will prove that
$g_{k}\dotsc g_{n}\cdot B^+$ lies in $\mathcal B_{\ge 0}$ for all
$k$ by descending induction. Suppose we know $g_{k+1}\dotsc
g_{n}\cdot B^+\in\mathcal B_{\ge 0}$. There are two possibilities
for $g_{k}$. The first case, $k\in J^\circ_{\v_+}$, is clear. In
that case $g_k=y_{i_k}(t_k)\in U^-_{\ge 0}$ and so $g_k
g_{k+1}\dotsc g_{n}\cdot B^+$ again lies in $\mathcal B_{\ge 0}$.

Let us now consider the other case. So $g_{k}=\dot s_{i_k}$ and
$k\in J^+_{\v_+}$. From \eqref{e:siIdentity} we get the formula
\begin{equation*}
x_i(t) \dot s_i=\alpha_i^\vee(t)y_i(t)x_i(-t\inv).
\end{equation*}
We apply this element for $i=i_k$ to $g_{k+1}\dotsc g_n\cdot B^+$
to get
 $$
x_{i_k}(t) \dot s_{i_k} g_{k+1}\dotsc g_n\cdot
B^+=\alpha_{i_k}^\vee(t)y_{i_k}(t) x_{i_k}(-t\inv) g_{k+1}\dotsc
g_n\cdot B^+.
 $$
Let $v'=v_{(k)}\inv v$ and $w'=w_{(k)}\inv w$. Since $\v$ is the
positive subexpression of $\w$, so the right-most reduced
subexpression and $k\in J^+_{\v}$, it follows that $s_{i_k}u>u$
for all $v'\le u\le w'$. Applying Lemma~\ref{l:tool} we get
\begin{equation*}
x_{i_k}(-t\inv) g_{k+1}\dotsc g_n\cdot B^+= g_{k+1}\dotsc g_n\cdot
B^+.
\end{equation*}
Therefore in total
\begin{equation}\label{e:finaltool}
 x_{i_k}(t) \dot s_{i_k} g_{k+1}\dotsc g_n\cdot B^+=
\alpha_{i_k}^\vee(t)y_{i_k}(t) g_{k+1}\dotsc g_n\cdot B^+.
\end{equation}
Now one can see that the right hand side of \eqref{e:finaltool} lies
in $\mathcal B_{\ge 0}$ for all $t>0 $. However as $t\to 0$ the
left hand side converges to $\dot s_{i_k} g_{k+1}\dotsc g_n\cdot
B^+$. Therefore $\dot s_{i_k}g_{k+1}\dotsc g_n\cdot B^+\in\mathcal
B_{\ge 0}$ and the lemma follows.
\end{proof}

\begin{proof}[Proof of Theorem~\ref{t:totpos}]
By Lemma~\ref{l:welldef} the map
\begin{equation}\label{e:map}
 G_{\v_+,\w}^{>0}\to\mathcal R^{>0}_{v,w}
\end{equation}
is well defined. Lemma~\ref{l:TotPosThenPos} implies that
$\mathcal R^{>0}_{v,w}\subset\mathcal R_{\v_+,\w}$. Now if
$B=z\dot w\cdot B^+\in\mathcal R^{>0}_{v,w}$, then
Lemma~\ref{l:oldlem} and Lemma~\ref{l:posminors} together with
Corollary~\ref{c:DeoInequalities} imply that the chamber
minors of $z$ (for the subexpression $\v_+$) are all positive.
Note also that $J^-_{\v_+}=\emptyset$. It follows that the map
$\mathcal R_{\v_+,\w}\to G_{\v_+,\w}$ described in Theorem~\ref{t:ansatz}
(see Proposition~5.2) restricts to
 $$
 \mathcal R^{>0}_{v,w}\to G_{\v_+,\w}^{>0},
 $$
giving the inverse to \eqref{e:map}.
\end{proof}

\section{Total positivity criteria}\label{s:PosCrit}
We can use the Chamber Ansatz together with Theorem~\ref{t:totpos}
to characterize $\mathcal R_{v,w}^{>0}$ by inequalities. In the
case where $v=1$ the criteria below reduce to the total positivity
criteria for $U^+\cap B^-\dot w B^-$ of Berenstein and Zelevinsky,
\cite{BeZel:TotPos}~Theorem~6.9.

\begin{prop}\label{p:PosCrit}
Consider $w\in W$ with a fixed reduced expression
$\w=(w_{(0)},\dotsc,w_{(n)})$. Let $v\le w$ in $W$ and
$\v_+=(v_{(0)},\dotsc, v_{(n)})$ be the positive subexpression
for $v$ in $\w$. Then
 \begin{align*}
 \mathcal R_{v,w}^{>0} &=\left \{z\dot w\cdot B^+\ \left |\ z \in
 U^+;
 \begin{array}{ll}\Delta^{v_{(k-1)}\omega_{i_k}}_{w_{(k)}\omega_{i_k}}(z)=0,&
 k\in J^+_{\v_+} \\
 \Delta^{v_{(k)}\omega_{i_k}}_{w_{(k)}\omega_{i_k}}(z)>0, & k\in
 J^\circ_{\v_+}
 \end{array}\right .\right \}\\&=
 \left\{ z\dot w\cdot B^+\in\mathcal R_{v,w}\left |\
 z\in U^+;\, \Delta^{v_{(k)}\omega_{i_k}}_{w_{(k)}\omega_{i_k}}(z)>0,
 \ k\in J^\circ_{\v_+}\right.\right\}.
 \end{align*}
\end{prop}

\begin{proof}
Let us call the two sets in question $S_1$ and $S_2$, respectively.
We want to show
$\mathcal R_{v,w}^{>0}=S_1=S_2$. The inclusions in one direction,
$R_{v,w}^{>0}\subseteq S_1\subseteq S_2$,
are clear from Remark~\ref{r:posminors} and Corollary~\ref{c:DeoInequalities}.

Moreover by Corollary~\ref{c:DeoInequalities},
\begin{equation*}
\Delta^{v_{(k-1)}\omega_{i_k}}_{w_{(k)}\omega_{i_k}}(z)=0, \quad
 \text{ for $k\in J^+_{\v_+}$,}
\end{equation*}
is true for all $z\dot w\cdot B^+$ in $\mathcal R_{\v_+,\w}$.
And since $\mathcal R_{\v_+,\w}$ is dense in $\mathcal R_{v,w}$,
this equality holds for all $z\dot w\cdot B^+\in\mathcal R_{v,w}$.
Therefore we also have the inclusions
$S_2\subseteq S_1\subseteq \mathcal R_{\v_+,\w}$.

It remains to show the final inclusion,
of $S_1$, say, into $\mathcal R_{v,w}^{>0}$.
Consider $z\dot w\cdot B^+\in\mathcal R_{\v_+,\w}$ with
$\Delta^{v_{(k)}\omega_{i_k}}_{w_{(k)}\omega_{i_k}}(z)>0$ for all
$k\in J^\circ_{\v_+}$. By Theorem~\ref{t:totpos} and
Theorem~\ref{t:ansatz}.(1) we need only show that the remaining
chamber minors,
$\Delta^{v_{(k)}\omega_{i_k}}_{w_{(k)}\omega_{i_k}}(z)$ for $k\in
J^+_{\v_+}$, are also positive. Suppose indirectly that
$\Delta^{v_{(k_0)}\omega_{i_{k_0}}}_{w_{(k_0)}\omega_{i_{k_0}}}(z)<0$
for some $k_0\in J^+_{\v_+}$. We may choose $k_0$ to be minimal with
this property. From Remark~\ref{r:formulas}.(1) and equation
\eqref{e:generalansatz} one obtains
 \begin{equation}\label{e:otherminors}
\Delta_{w_{(k_0)}\omega_{i_{k_0}}}^{v_{(k_0)}\omega_{i_{k_0}}}(z)
=\frac{
\prod_{j\not=i_{k_0}}
\Delta_{w_{(k_0)}\omega_j}^{v_{(k_0)}\omega_j}(z)^{-a_{j,i_{k_0}}}}
{\Delta_{w_{(k_0-1)}\omega_{i_{k_0}}}^{v_{({k_0}-1)}\omega_{i_{k_0}}}(z)}.
 \end{equation}
Since the right hand side is made up of chamber minors for
smaller $k$, it must be positive. So we have a contradiction.
\end{proof}

\begin{rem}
Note that by this proposition, $\mathcal R_{v,w}^{>0}$ is given
inside $\mathcal R_{v,w}$ by $\dim (\mathcal
R_{v,w})=\ell(w)-\ell(v)$ inequalities. This is the ideal
situation. Indeed, it is easy to see that our set of inequalities,
 \begin{equation}\label{e:ineq}
 \Delta^{v_{(k)}\omega_{i_{k}}}_{w_{(k)}\omega_{i_k}}(z)>0\qquad
 (k\in J^\circ_{\v_+}),
 \end{equation}
is minimal. Using Proposition~\ref{p:ansatz} we can describe the
situation by the following commutative diagram
\begin{equation*}
\begin{array}{cccl}
 \Delta:& \mathcal R_{\v_+,\w} &\overset{\sim}\To & (\R^*)^{\ell(w)-\ell(v)} \\
        &   \uparrow   &                     &\quad \uparrow \\
        & \mathcal R_{v,w}^{>0} &\overset{\sim}\To &
(\R_{>0})^{\ell(w)-\ell(v)},
\end{array}
\end{equation*}
where the horizontal maps are given by the chamber minors,
 $$\Delta\left (z\dot w\cdot B^+\right)=
\left (\Delta^{v_{(k)}\omega_{i_{k}}}_{w_{(k)}\omega_{i_k}}(z)
\right)_{k\in J^\circ_{\v_+}},$$
and the vertical maps are inclusions.
>From this picture it is clear that \eqref{e:ineq}, or
Proposition~\ref{p:PosCrit},
has no redundant inequalities.
\end{rem}

\noindent {\bf Acknowledgements.} \\
The first author would like to thank King's College, London for their
kind hospitality in spring 2003, where most of the work for
this paper was carried out.



\bibliographystyle{plain}
\def\cprime{$'$}

\end{document}